\documentstyle[amscd,amssymb,verbatim,12pt]{amsart}
\newcommand{\aff}{\mbox{\rm aff}}
\newcommand{\Aff}{\mbox{\rm Aff}}

\newcommand{\Aut}{\mbox{\rm Aut}}
\newcommand{\bb}{\mbox{\rm b}}

\newcommand{\Coker}{\mbox{\rm Coker}}

\newcommand{\const}{\mbox{\rm const.}}

\newcommand{\diam}{\mbox{\rm diam}}
\newcommand{\Diff}{\mbox{\rm Diff}}

\newcommand{\End}{\mbox{\rm End}}

\newcommand{\GL}{\mbox{\rm GL}}

\newcommand{\HH}{\mbox{\rm H}}

\newcommand{\Hom}{\mbox{\rm Hom}}

\newcommand{\Id}{\mbox{\rm Id}}

\newcommand{\Image}{\mbox{\rm Im}}

\newcommand{\inj}{\mbox{\rm inj}}

\newcommand{\Ker}{\mbox{\rm Ker}}

\newcommand{\N}{{\Bbb N}}

\newcommand{\Out}{\mbox{\rm Out}}

\newcommand{\R}{{\Bbb R}}

\newcommand{\Ric}{\mbox{\rm Ric}}

\newcommand{\vol}{\mbox{\rm vol}}
\newcommand{\Z}{{\Bbb Z}}
\theoremstyle{plain}
\newtheorem{definition}{Definition}

\newtheorem{lemma}{Lemma}
\newtheorem{theorem}{Theorem}
\newtheorem{proposition}{Proposition}
\newtheorem{corollary}{Corollary}

\numberwithin{equation}{section}
\renewcommand{\rm}{\normalshape}
\begin{document}
\title{Collapsing and the Differential Form Laplacian : The Case of a 
Smooth Limit Space}
\author{John Lott}
\address{Department of Mathematics\\
University of Michigan\\
Ann Arbor, MI  48109-1109\\
USA}
\email{lott@@math.lsa.umich.edu}
\thanks{Research supported by NSF grant DMS-9704633}
\subjclass{Primary: 58G25; Secondary: 53C23}
\date{May 20, 2000}
\maketitle
\begin{abstract}
We analyze the limit of the $p$-form
Laplacian under a collapse, with bounded sectional curvature and bounded
diameter, to a smooth limit space.  
As an application, we
characterize when the $p$-form 
Laplacian has small positive eigenvalues in a collapsing sequence.
\end{abstract}

\section{Introduction} \label{sect1}

A central problem in geometric analysis is to estimate the spectrum of the
Laplacian on a compact Riemannian
manifold $M$ in terms of geometric invariants.  In the case of
the Laplacian on functions, a major result is Cheeger's lower bound on the
smallest positive eigenvalue in terms of an isoperimetric constant
\cite{Cheeger (1970)}. The problem of extending his lower bound
to the case of the $p$-form Laplacian was posed in \cite{Cheeger (1970)}.  
There has been little
progress on this problem.  We will address the more general question of 
estimating the eigenvalues $\{\lambda_{p,j}(M)\}_{j=1}^\infty$
of the $p$-form Laplacian $\triangle_p$ (counted with multiplicity) 
in terms of geometric invariants of $M$.

A basic fact, due to Cheeger and Dodziuk, 
is that $\lambda_{p,j}(M)$ depends continuously on the
Riemannian metric $g^{TM}$ in the $C^0$-topology \cite{Dodziuk (1982)}.
Then an immediate consequence of the $C^\alpha$-compactness theorem of
Anderson and Cheeger \cite{Anderson-Cheeger (1992)} is that for any
$n \in \Z^+$, $r \in \R$, and $D, i_0 > 0$, 
there are uniform 
bounds on $\lambda_{p,j}(M)$ among connected closed $n$-dimensional
Riemannian manifolds $M$ with $\Ric(M) \: \ge \: r$, $\diam(M) \le D$
and $\inj(M) > i_0$
(compare \cite[Theorem 1.3]{Chanillo-Treves (1997)},
\cite[Theorem 0.4]{Colbois-Courtois (1990)}.) In particular, there is
a uniform positive lower bound on the smallest positive
eigenvalue of the $p$-form
Laplacian under these geometric assumptions. 

The question, then, is what happens when
$\inj(M) \rightarrow 0$.
For technical reasons, in this paper we will
assume uniform bounds on the Riemannian curvature $R^M$. Then we
wish to study how the spectrum of
$\triangle_p$ behaves in the collapsing limit.
By collapsing we mean the phenomenon of
a sequence of Riemannian manifolds converging in the Gromov-Hausdorff
topology to a lower-dimensional space. We refer to
\cite[Chapters 1 and 3]{Gromov (1999)} for basic information about
collapsing and 
\cite[Section I]{Cheeger-Fukaya-Gromov (1992)},
\cite{Cheeger-Gromov (1986)},
\cite{Fukaya (1990)} and
\cite[Chapter 6]{Gromov (1999)} for information about bounded
curvature collapsing.
In this paper, we analyze the behavior of the spectrum of $\triangle_p$ under
collapse, with bounded sectional curvature and bounded diameter, to a
smooth limit space. The answer
will be in terms of a type of Laplacian on the limit space.  
As an application, we 
characterize when the $p$-form 
Laplacian has small positive eigenvalues in a collapsing sequence.
In a subsequent paper we will extend the results to the case of singular
limit space, and give additional applications.

From Hodge theory, $\dim(\Ker(\triangle_p)) = \bb_p(M)$, the $p$-th
Betti number of $M$. Given $K \: \ge \: 0$, let ${\cal M}(M, K)$ be the
set of Riemannian metrics $g$ on $M$ with $\parallel R^{M} \parallel_\infty
\: \le \: K$ and $\diam(M, g) \: \le \: 1$.
We will say that $M$ has small positive
eigenvalues of the $p$-form Laplacian if  
\begin{equation}
\inf_{g \in {\cal M}(M, K)} \lambda_{p,j}(M,g) \: = \: 0
\end{equation}
for some $j \: > \: \bb_p(M)$ and
some $K > 0$.
If this is the case then we will say
that $M$ has (at least) $j$ small eigenvalues. Note that this is a statement
about the (smooth) topological type of $M$.

There are no small positive eigenvalues of the
Laplacian on functions on $M$ (see, for example,
\cite{Berard (1988)}). 
Colbois and Courtois gave examples of manifolds with
small positive eigenvalues of the $p$-form Laplacian for $p > 0$
\cite{Colbois-Courtois (1990)}.
Their examples were manifolds $M$ with free isometric
$T^k$-actions, which one shrinks
in the direction of the $T^k$-orbits. In terms of the
fiber bundle  $M \rightarrow M/T^k$, this sort of
collapsing is a case of the so-called adiabatic limit. The asymptotic
behaviour of the small
eigenvalues of the $p$-form Laplacian in the adiabatic limit
was related to the Leray spectral sequence of the fiber
bundle in 
\cite{Berthomieu-Bismut (1994),Dai (1991),Forman (1995),Mazzeo-Melrose (1990)}.

In another direction, Fukaya considered the behavior of the Laplacian on
functions in the case of a sequence of manifolds that converge in the
Gromov-Hausdorff metric $d_{GH}$ to a lower-dimensional
limit space $X$, the collapsing assumed to be with
bounded sectional curvature and bounded diameter \cite{Fukaya (1987)}.
He found that in order to get limits, one needs to widen the class of spaces 
being considered
by adding a Borel measure, and consider measured metric spaces.
This is the case even if $X$ happens to be a smooth manifold. 
He defined a Laplacian acting
on functions on the measured limit space and proved a convergence theorem 
for the spectrum of the Laplacian on functions, under the geometric
assumption of convergence in the measured Gromov-Hausdorff topology.

We consider the behavior of the spectrum of $\triangle_p$ under
collapse with bounded sectional curvature and bounded diameter. 
We find that we need a somewhat
more refined structure on the limit space, namely a superconnection as
introduced by Quillen \cite{Quillen (1985)}. 
More precisely, we will need a flat
degree-$1$ superconnection in the sense of \cite{Bismut-Lott (1995)}.
Suppose that
$B$ is a smooth connected closed manifold and that $E = \bigoplus_{j=0}^m
E^j$ is a $\Z$-graded real vector bundle on $B$. The degree-$1$
superconnections $A^\prime$ that we need will be of the form
\begin{equation} \label{eq1.3}
A^\prime \: = \: A^\prime_{[0]} \: +  \: A^\prime_{[1]} \: + \: 
A^\prime_{[2]}
\end{equation}
where
\begin{itemize}
\item $A^\prime_{[0]} \in C^\infty \left( B; \Hom(E^*, E^{*+1}) \right)$,
\item $A^\prime_{[1]}$ is a grading-preserving connection $\nabla^E$ on $E$ 
and 
\item $A^\prime_{[2]} \in \Omega^2 \left( B; \Hom(E^*, E^{*-1}) \right)$.
\end{itemize}
The superconnection extends by Leibniz' rule to an operator $A^\prime$
on the $E$-valued differential forms
$\Omega (B; E)$.
The flatness condition $\left( A^\prime \right)^2 = 0$ becomes
\begin{itemize}
\item
$\left( A^\prime_{[0]} \right)^2 = \left( A^\prime_{[2]} \right)^2 = 0$,
\item $\nabla^E A^\prime_{[0]} = \nabla^E A^\prime_{[2]} = 0$ and
\item $\left( \nabla^E \right)^2 + A^\prime_{[0]} A^\prime_{[2]} +
A^\prime_{[2]} A^\prime_{[0]} = 0$.
\end{itemize}
In particular, $A^\prime_{[0]}$ defines a differential complex on
the fibers of $E$.
Let $g^{TB}$ be a Riemannian metric on $B$ and let $h^E$ be a graded Euclidean
inner product on $E$, meaning that $E^i$ is orthogonal to $E^j$ if $i \ne j$.
Then there are an adjoint $\left( A^\prime \right)^*$ 
to $A^\prime$ and a Laplacian
$\triangle^E = A^\prime \left( A^\prime \right)^* + \left( A^\prime \right)^*
A^\prime$ on $\Omega (B; E)$. Let $\triangle^E_p$ be the restriction of
$\triangle^E$ to $\bigoplus_{a+b=p} \Omega^a(B; E^b)$.

Using the $C^0$-continuity of the
spectrum and the geometric results of Cheeger, Fukaya and Gromov
\cite{Cheeger-Fukaya-Gromov (1992)}, 
we can reduce our study of collapsing to certain special fiber
bundles. As is recalled in Section \ref{sect3},
 an infranilmanifold $Z$ has a canonical flat
linear connection $\nabla^{aff}$. Let $\Aff(Z)$ be the group of
diffeomorphisms of $Z$ which preserve $\nabla^{aff}$.

\begin{definition} \label{def2}
An affine fiber bundle is a smooth fiber bundle $M \rightarrow B$ whose
fiber $Z$ is an infranilmanifold and whose structure group is reduced
from $\Diff(Z)$ to $\Aff(Z)$.
A Riemannian affine fiber bundle is an affine fiber bundle along with
\begin{itemize}
\item A horizontal distribution $T^HM$ on $M$ whose holonomy lies in 
$\Aff(Z)$,
\item
A family $g^{TZ}$ of vertical Riemannian metrics which
are parallel with respect to the flat affine connections on the 
fibers $Z_b$ and
\item
A Riemannian metric $g^{TB}$ on $B$.
\end{itemize}
\end{definition}

Fix a smooth connected closed Riemannian manifold $B$.
Fukaya showed that any manifold $M$ which collapses to
$B$, with bounded sectional curvature, is the
total space of an affine fiber bundle over $B$ \cite{Fukaya (1989)}.
If $M \rightarrow B$ is an affine fiber bundle, let $T^HM$ be a horizontal
distribution on $M$ as above.
Let $T \in \Omega^2(M; TZ)$ be
the curvature of $T^HM$. There is a
$\Z$-graded real vector bundle $E$ on $B$ whose fiber over $b \in B$ is
isomorphic to the differential forms on the fiber $Z_b$ which are
parallel with respect to the flat affine connection on
$Z_b$. The exterior derivative $d^M$ induces a
flat degree-$1$ superconnection $A^\prime$ on $E$.
If $M \rightarrow B$ is in addition a Riemannian affine fiber bundle then 
we obtain a Riemannian metric $g^{TM}$ on $M$ 
constructed from $g^{TZ}$, $g^{TB}$ and $T^HM$. 
There is 
an induced $L^2$-inner product $h^E$ on $E$.
Define $\triangle^E$ as above.
Let $\diam(Z)$ denote the maximum diameter of the fibers 
$\{Z_b\}_{b \in B}$ in the intrinsic metric
and let $\Pi$ denote the second fundamental
forms of the fibers $\{Z_b\}_{b \in B}$.
Our first result says 
that  the spectrum $\sigma(\triangle^E_p)$ of $\triangle^E_p$
contains all of the spectrum of the 
$p$-form Laplacian $\triangle^M_p$ which stays bounded
as $d_{GH}(M, B) \rightarrow 0$.

\begin{theorem} \label{th1}
There are positive constants $A$, $A^\prime$ and 
$C$ which only depend on $\dim(M)$ such that
if $\parallel R^Z \parallel_\infty \diam(Z)^2 \le A^\prime$ then
for all $0 \: \le \: p \: \le \dim(M)$,
\begin{align} \label{eq1.4}
& \sigma(\triangle^M_p) \cap \left[ 0, \: A \: 
{\diam(Z)^{-2}} \: - C \: \left(
\parallel R^{M} \parallel_\infty + 
\parallel \Pi \parallel_\infty^2 +
\parallel T \parallel_\infty^2 
 \right) \right) = \\
& \sigma(\triangle^E_p) \cap \left[ 0, \: A \:
{\diam(Z)^{-2}} \: - C \: \left(
\parallel R^{M} \parallel_\infty + 
\parallel \Pi \parallel_\infty^2 + 
\parallel T \parallel_\infty^2 
\right) \right). \notag
\end{align}
\end{theorem}

When $Z$ is flat, there is some intersection between Theorem \ref{th1} and
the adiabatic limit results of 
\cite{Berthomieu-Bismut (1994),Dai (1991),Forman (1995),Mazzeo-Melrose (1990)}.
However, there is
the important difference that we need estimates which are uniform
with respect to $d_{GH}(M, B)$, whereas the adiabatic limit results concern
the asymptotics of the eigenvalues under the collapse of a given 
Riemannian fiber bundle coming from a constant rescaling of its fibers.

We apply Theorem \ref{th1} to estimate the eigenvalues of a 
general Riemannian manifold $M$
which is Gromov-Hausdorff close to $B$, assuming sectional curvature bounds
on $M$.
Of course we cannot say precisely what $\sigma(\triangle^M)$ is, but we
can use Theorem \ref{th1} 
to approximate it to a given precision $\epsilon > 0$.
We say that two nonnegative numbers $\lambda_1$ and $\lambda_2$ are
$\epsilon$-close if
$e^{-\epsilon} \: \lambda_{2} \: \le \: \lambda_{1} \: \le \: e^\epsilon \: 
\lambda_{2}$.
We show that for a given $\epsilon > 0$, if $d_{GH}(M, B)$
is sufficiently small
then there is a flat degree-$1$ superconnection $A^\prime$ on $B$ whose 
Laplacian $\triangle^E_p$ has a spectrum which is $\epsilon$-close to that of
$\triangle^M_p$, at least up to a high level.

\begin{theorem} \label{th2}
Let $B$ be a fixed
smooth connected closed
Riemannian manifold.  Given $n \in \Z^+$, $\epsilon > 0$ and
$K \ge 0$, there are positive constants $A(n, \epsilon, K)$,
$A^\prime(n, \epsilon, K)$ and 
$C(n, \epsilon, K)$ with the following property :
If $M^n$ is an $n$-dimensional connected closed
Riemannian manifold with $\parallel R^M \parallel_\infty
\: \le \: K$ and $d_{GH}(M, B) \: \le 
\: A^\prime(n,\epsilon, K)$
then there are\\
1. A $\Z$-graded real vector bundle $E$ on $B$,\\
2. A flat degree-$1$ superconnection $A^\prime$ on $E$ and\\
3. A Euclidean inner product $h^E$ on $E$\\
such that    
if $\lambda_{p,j}(M)$ is the $j$-th eigenvalue of the $p$-form Laplacian
on $M$, $\lambda_{p,j}(B; E)$ is the $j$-th eigenvalue of
$\triangle^E_p$ and
\begin{equation} \label{eq1.5}
\min (\lambda_{p,j}(M), \lambda_{p,j}(B; E)) \: \le \: A(n,\epsilon, K) \: 
{d_{GH}(M, B)^{-2}} \: - C(n,\epsilon, K)
\end{equation}
then $\lambda_{p,j}(M)$ 
is $\epsilon$-close to $\lambda_{p,j}(B; E)$.
\end{theorem}

Using \cite{Baker-Dodziuk (1997)}, one can also show that
the eigenspaces of $\triangle^E_p$ are $L^\infty$-close to those of 
$\triangle^M_p$, with respect to the embedding $\Omega (B; E) \rightarrow
\Omega(M)$.

In the case of the Laplacian on functions, only $E^0$ is relevant.
Although $E^0$ is the trivial $\R$-bundle on $B$ with a trivial connection, 
its Euclidean inner product
$h^{E^0}$ need not be trivial and corresponds exactly to the measure
in Fukaya's work. 

In order to apply Theorem \ref{th2}, we prove a compactness result for the
superconnection and Euclidean metric.

\begin{definition} \label{def3}
Let ${\cal S}_E$ be the space of degree-$1$
superconnections on $E$,
let ${\cal G}_E$ be the group of smooth grading-preserving
$\GL(E)$-gauge transformations on $E$ and
let ${\cal H}_E$ be the space of graded Euclidean inner products on 
$E$. We equip ${\cal S}_E$ and ${\cal H}_E$
with the $C^\infty$-topology. Give $({\cal S}_{E} \times {\cal H}_{E})/
{\cal G}_{E}$ the quotient topology. 
\end{definition}

\begin{theorem} \label{th3}
In Theorem \ref{th2}, we may assume that $E$ is one of a finite number
of isomorphism classes of
real $\Z$-graded topological vector bundles $\{E_i\}$ on $B$. 
Furthermore, there are compact subsets $D_{E_i} \subset
({\cal S}_{E_i} \times {\cal H}_{E_i})/{\cal G}_{E_i}$ 
depending on $n$, $\epsilon$ and $K$, such that we 
may assume that the gauge-equivalence class of the pair
 $\left(A^\prime,h^E \right)$
lies in $D_E$.
\end{theorem}

We remark that there may well 
be a sequence of topologically distinct Riemannian manifolds of a
given dimension,
with uniformly bounded sectional curvatures, which converge to $B$ in the
Gromov-Hausdorff topology (see Example 3 of Section \ref{sect2}).
This contrasts with the finiteness statement in Theorem \ref{th3}.

The eigenvalues of $\triangle^E_p$ are continuous with respect to 
$\left[ (A^\prime,h^E \right)] \in ({\cal S}_{E} \times {\cal H}_{E})/
{\cal G}_{E}$.
One application of Theorem \ref{th3} is the following relationship between
the spectra of $\triangle^M_p$ and the ordinary differential form
Laplacian on $B$.

\begin{theorem} \label{th4}
Under the hypotheses of
Theorem \ref{th2}, let $\lambda^\prime_{p,j}(B)$ be the
$j$-th eigenvalue of the Laplacian on $\bigoplus_r \Omega^r(B) \otimes
\R^{dim(E^{p-r})}$.  Then there is a positive 
constant $D(n, \epsilon, K)$ such that
\begin{equation} \label{eq1.6}
e^{-\epsilon/2} \: \lambda^\prime_{p,j}(B)^{1/2} \: - \: D(n, \epsilon, K) \:
\le \lambda_{p,j}(M)^{1/2} \: \le \: 
e^{\epsilon/2} \: \lambda^\prime_{p,j}(B)^{1/2} \: + \: D(n, \epsilon, K).
\end{equation}
\end{theorem}

Now consider a flat degree-$1$
superconnection $A^\prime$ on a real $\Z$-graded
vector bundle $E$ over a smooth manifold $B$.
As $\left( A^\prime \right)^2 = 0$, there is a
cohomology $\HH^*(A^\prime)$ for the action of $A^\prime$ on 
$\Omega (B; E)$, the latter having the total grading. There is a flat
$\Z$-graded
``cohomology'' vector bundle $\HH^*(A^\prime_{[0]})$ on $B$. 
Furthermore, there is a spectral sequence to compute $\HH^*(A^\prime)$, with
$E_2$-term $\HH^* \left(B; \HH^*(A^\prime_{[0]}) \right)$.

Suppose that $M$ is a connected closed manifold with at least $j$ small
eigenvalues of $\triangle_p$ for
$j > \bb_p(M)$. Consider a
sequence of Riemannian metrics $\{g_i\}_{i=1}^\infty$ in ${\cal M}(M, K)$ 
with $\lim_{i \rightarrow\infty} 
\lambda_{p,j}(M, g_i) = 0$. There must be a subsequence of
$\{(M, g_i)\}_{i=1}^\infty$
which converges to a lower-dimensional limit space $X$. That is, we are in 
the collapsing situation. 
Suppose that the limit space is
a smooth manifold $B$.
From Theorems \ref{th2} and \ref{th3},
we can take a further subsequence of $\{(M, g_i)\}_{i=1}^\infty$ 
to obtain a single vector bundle $E$ on $B$, equipped
with a sequence $\left\{ \left( A^\prime_i,
h^E_i \right) \right\}_{i=1}^\infty$ of superconnections and Euclidean inner
products.
Using the compactness result in Theorem \ref{th3}, we can take a convergent
subsequence of these pairs, modulo gauge transformations, to obtain a 
superconnection $A^\prime_\infty$ on $E$ with 
$\dim \Ker \left( \triangle^E_p \right) \: \ge \: j$.
Then $\dim \left( \HH^p(A^\prime_\infty) \right) \: \ge \: j$.
It is no longer true that 
$\HH^*(A^\prime_\infty) \cong \HH^*(M ; \R)$ for this limit superconnection.
However, we can analyze $\HH^*(A^\prime_\infty)$ using the spectral sequence.
We obtain
\begin{equation} \label{eq1.7}
j \: \le \: \sum_{a+b=p} \dim \left( \HH^a(B; \HH^b(A_{\infty,[0]}^\prime)) 
\right).
\end{equation}

This formula has some immediate consequences.
The first one is a bound on the number of
small eigenvalues of the $1$-form Laplacian.

\begin{corollary} \label{cor2} 
Suppose that $M$ has $j$ small eigenvalues of the $1$-form Laplacian,
with $j \: > \: \bb_1(M)$. Let $X$ be the limit space
coming from the above argument. Suppose that $X$ is a smooth manifold $B$.
Then
\begin{equation}
j \: \le \:
\bb_1(B) \: + \: \dim(M) \: - \: \dim(B) \: \le \: \bb_1(M) \: + \: \dim(M).
\end{equation} 
\end{corollary}

The second consequence
is a bound on the number of small eigenvalues of the $p$-form
Laplacian for a manifold which is Gromov-Hausdorff close to a codimension-$1$
manifold.

\begin{corollary} \label{cor4}
Let $B$ be a connected closed $(n-1)$-dimensional Riemannian
manifold.  Then for any $K \ge 0$, there are
$\delta, c > 0$ with the following property :
Suppose that $M$ is a connected closed smooth $n$-dimensional 
Riemannian manifold with $\parallel R^M \parallel_\infty \: \le \: K$
and $d_{GH}(M, B) < \delta$. First, $M$ is the total space of
a circle bundle over $B$.  Let ${\cal O}$ be the
orientation bundle of $M \rightarrow B$, a flat real line bundle on
$B$. Then
$\lambda_{p,j}(M, g) \: > \: c$ for 
$j \: = \: \bb_p(B) + \bb_{p-1}(B; {\cal O}) \: + \: 1$.
\end{corollary}

The rest of our results concern small eigenvalues in collapsing sequences.

\begin{definition}
If $M \rightarrow B$ is an affine fiber bundle, 
a collapsing sequence associated to the affine fiber bundle
is a sequence of metrics $\{g_i\}_{i=1}^\infty \in
{\cal M}(M, K)$ for some $K \ge 0$
such that $\lim_{i \rightarrow \infty} (M, g_i) = B$ in the 
Gromov-Hausdorff topology and for some $\epsilon > 0$, each $(M, g_i)$ is 
$\epsilon$-biLipschitz to a Riemannian affine fiber bundle structure on
$M \rightarrow B$.
\end{definition}

We show that there are three mechanisms to make small positive eigenvalues
of the differential form Laplacian on $M$ in a collapsing sequence. 
Either the differential form Laplacian on the fiber admits small positive
eigenvalues, or
the holonomy of the flat ``cohomology'' bundle on $B$ fails to
be semisimple, or the Leray spectral sequence of $M \rightarrow B$ does not
degenerate at the $E_2$-term.

\begin{theorem} \label{th6}
Let $\{(M, g_i)\}_{i=1}^\infty$ 
be a collapsing sequence associated to an affine fiber
bundle $M \rightarrow B$. Suppose that 
$\lim_{i \rightarrow \infty} \lambda_{p,j}(M, g_i) = 0$ for
some $j > \bb_p(M)$. 
Write the fiber $Z$ of the affine fiber bundle as the quotient of a
nilmanifold $\widehat{Z} \: = \: \widehat{\Gamma} \backslash N$ by a finite
group $F$.
Then \\
1. For some $q \in [0,p]$, $b_q(Z) \: < \: \dim \left(\Lambda^q({\frak n}^*)^F 
\right)$, or\\
2. For all $q \in [0,p]$, $b_q(Z) \: = \: \dim \left(\Lambda^q({\frak n}^*)^F 
\right)$, and for 
some $q \in [0, p]$, the holonomy representation of the flat vector
bundle $\HH^q(Z; \R)$ on $B$ fails to be semisimple, or\\
3. For all $q \in [0,p]$, $b_q(Z) \: = \: \dim \left(\Lambda^q({\frak n}^*)^F 
\right)$ and the holonomy representation of the flat vector
bundle $\HH^q(Z; \R)$ on $B$ is semisimple, and
the Leray spectral sequence to compute $\HH^p(M; \R)$
does not degenerate at the $E_2$ term.
\end{theorem}

Examples show that
small positive
eigenvalues can occur in each of the three cases in Theorem \ref{th6}.

Theorem \ref{th6} has some immediate consequences.  The first is 
a characterization of when the $1$-form Laplacian has small
positive eigenvalues in
a collapsing sequence.

\begin{corollary} \label{cor5}
Let $\{(M, g_i)\}_{i=1}^\infty$ 
be a collapsing sequence associated to an affine fiber
bundle $M \rightarrow B$. Suppose that 
$\lim_{i \rightarrow \infty} \lambda_{1,j}(M, g_i) = 0$ for
some $j > \bb_1(M)$. Then \\
1. The differential $d_2 : 
\HH^0(B; \HH^1(Z; \R)) \rightarrow \HH^2(B; \R)$ in the Leray spectral sequence
for $\HH^*(M; \R)$ is nonzero, or\\
2. The holonomy representation of the flat vector
bundle $\HH^1(Z; \R)$ on $B$ has a nontrivial unipotent subrepresentation, or\\
3. $Z$ is almost flat but not flat and there is a nonzero covariantly-constant
section of the flat vector bundle 
$\frac{H^1(A^\prime_{\infty,[0]})}{H^1(Z; \R)_\infty}$.
\end{corollary}

The differential $d_2 : 
\HH^0(B; \HH^1(Z; \R)) \rightarrow \HH^2(B; \R)$ can be considered to be a
type of Euler class; in the case of an oriented circle bundle over a smooth
base, it gives exactly the Euler class.

The second consequence is
a characterization of when the $p$-form Laplacian has small 
positive eigenvalues in
a collapsing sequence over a circle.

\begin{corollary} \label{cor6}
Let $\{(M, g_i)\}_{i=1}^\infty$ 
be a collapsing sequence associated to an affine fiber
bundle $M \rightarrow S^1$.
Suppose that 
$\lim_{i \rightarrow \infty} \lambda_{p,j}(M, g_i) = 0$ for
some $j > \bb_p(M)$.  Write the fiber $Z$ of the affine fibre bundle
as in Theorem \ref{th6}. Then \\
1. For some $q \in \{p-1, p\}$, $b_q(Z) \: < \: \dim \left(\Lambda^q({\frak n}^*)^F 
\right)$, or\\
2. For $q \in \{p-1, p\}$, $b_q(Z) \: = \: \dim \left(\Lambda^q({\frak n}^*)^F 
\right)$, and if $\Phi^* \in \Aut(\HH^*(Z; \R))$ denotes the holonomy
action on the fiber cohomology then $\Phi^p$ or $\Phi^{p-1}$ has a
nontrivial unipotent factor in its Jordan normal form.
\end{corollary}

The third consequence
is a characterization of when the $p$-form Laplacian has small
positive eigenvalues in a collapsing sequence over a codimension-$1$ 
manifold.

\begin{corollary} \label{cor7}
Let $\{(M, g_i)\}_{i=1}^\infty$ 
be a collapsing sequence associated to an affine fiber
bundle $M \rightarrow B$ with $\dim(B) = \dim(M) - 1$. Suppose that 
$\lim_{i \rightarrow \infty} \lambda_{p,j}(M, g_i) = 0$ for
some $j > \bb_p(M)$.  Let ${\cal O}$ be the
orientation bundle of $M \rightarrow B$, a flat real line bundle on
$B$.
Let $\chi \in \HH^2(B; {\cal O})$ be the Euler class of the
orbifold circle bundle $M \rightarrow B$. Let ${\cal M}_\chi$ be multiplication
by $\chi$. Then
${\cal M}_\chi : 
\HH^{p-1}(B; {\cal O}) \rightarrow \HH^{p+1}(B; \R)$ is nonzero
or ${\cal M}_\chi : 
\HH^{p-2}(B; {\cal O}) \rightarrow \HH^{p}(B; \R)$ is nonzero.
\end{corollary}

Finally, we give a class of examples for which the inequality in
(\ref{eq1.7}) is an equality.

\begin{theorem} \label{th7}
Suppose that $M \rightarrow B$ is a affine fiber bundle with a smooth
base $B$ and fiber $Z = \widehat{Z}/F$, where $\widehat{Z}$ is a nilmanifold
$\widehat{\Gamma} \backslash N$ and $F$ is a finite group. 
Let
\begin{equation} \label{eq1.8}
{\frak n} = {\frak n}^\prime_{[0]} \supset {\frak n}^\prime_{[1]} \supset 
\ldots \supset {\frak n}^\prime_{[S]} \supset 0
\end{equation}
be the lower central series of the Lie algebra ${\frak n}$. Let 
${\frak c}({\frak n})$ be the center of ${\frak n}$. For $0 \le k \le S$, put
\begin{equation} \label{eq1.9}
{\frak n}_{[k]} \: = \: {\frak n}^\prime_{[k]} \: + \: {\frak c}({\frak n})
\end{equation}
and put ${\frak r}_{[k]} =
{\frak n}_{[k]}/{\frak n}_{[k+1]}$.
Let $P$ be the principal $\Aff(Z)$-bundle such that $M = P \times_{Aff(Z)} Z$. 
Let $G = \bigoplus_b G^b$ be the 
$\Z$-graded flat vector bundle on $B$ with 
\begin{equation} \label{eq1.10}
G^b = P \times_{Aff(Z)} \left(
\Lambda^b \left(  \bigoplus_{k=0}^S {\frak r}_{[k]}^* \right) \right)^F.
\end{equation}
Then for any $0 \le p \le \dim(M)$, 
$M$ has $\sum_{a+b=p} \dim
\left( \HH^a(B; G^b) \right) $ small eigenvalues of the $p$-form Laplacian.
\end{theorem}

The structure of the paper is as follows.  In Section \ref{sect2} we give
examples of collapsing which show that the superconnection formalism is
necessary.
In Section \ref{sect3} we give some background information about 
infranilmanifolds $Z$ and show that
the orthogonal
projection onto the parallel forms of $Z$ is independent of the choice of 
parallel
Riemannian metric.
In Section \ref{infranil} we give a detailed analysis of the spectrum of
the differential form Laplacian on an infranilmanifold.
In Section \ref{sect4} we show that the eigenvalues of the superconnection
Laplacian are continuous with respect to the superconnection, 
the Riemannian metric and the Euclidean inner product. 
We then analyze the differential form Laplacian on a Riemannian 
affine fiber bundle
and prove Theorem \ref{th1}.
In Section \ref{sect5} 
we consider manifolds $M$ that are Gromov-Hausdorff close to a smooth 
manifold $B$ and prove 
Theorems \ref{th2}, \ref{th3} and \ref{th4}. 
Section \ref{sect8} 
uses the compactness results to prove
Theorem \ref{th6} and Corollaries \ref{cor2}-\ref{cor7}. We then prove
Theorem \ref{th7}. More detailed descriptions appear at the
beginnings of the sections

After this paper was finished, I learned of the preprint version of
\cite{Colbois-Courtois (1998)} which, among other things, contains proofs of
Corollaries \ref{cor4} and \ref{cor7} 
in the case when $M$ and $B$
are oriented. The paper
\cite{Jammes (2000)} is also related to the present paper.

I thank Bruno Colbois, Gilles Courtois and Pierre Jammes for 
corrections to an earlier version of this paper. I thank the referee
for a very careful reading of the manuscript and many useful remarks,
along with suggesting a simplification of the proof of
Proposition \ref{prop2}.

\section{Examples} \label{sect2}

As for notation in this paper, 
if $G$ is a group which acts on a set $X$, we let $X^G$ denote the
set of fixed-points. If $B$ is a smooth manifold and $E$ is a smooth
vector bundle on $B$, we let $\Omega(B; E)$ denote the smooth
$E$-valued differential forms on $B$. If ${\frak n}$ is a nilpotent
Lie algebra on which a finite group $F$ acts by automorphisms
then ${\frak n}^*$ denotes the dual space, $\Lambda^*({\frak n}^*)$ denotes
the exterior algebra of the dual space and $\Lambda^*({\frak n}^*)^F$ denotes
the $F$-invariant subspace of the exterior algebra,  \\ \\
\noindent
{\bf Example 1 :} Let $N$ be a simply-connected connected nilpotent Lie
group, such as the $3$-dimensional Heisenberg group. Let ${\frak n}$ be
its Lie algebra of left-invariant vector fields, let
$g^{TN}$ be a left-invariant Riemannian metric on $N$ and let $\triangle^N$
be the corresponding Laplacian on $\Omega^*(N)$. (For simplicity of notation, 
we omit reference to the form degree $p$.)
The left-invariant differential forms
$\Lambda^*({\frak n}^*)$ form a subcomplex of $\Omega^*(N)$ with differential
$d^{\frak n}$, on which
$\triangle^N$ restricts to a finite-dimensional operator 
$\triangle^{\frak n}$. If $\Gamma$ is a lattice in $N$ then
the left-invariant forms on $N$ push down to forms on 
$Z = \Gamma \backslash N$, giving a subcomplex of $\Omega^*(Z)$ which is
isomorphic to $\Lambda^*({\frak n}^*)$. One knows that
$\HH^*(Z; \R)$ is isomorphic to the cohomology of this subcomplex
\cite[Corollary 7.28]{Raghunathan (1972)}. 
We see that the spectrum $\sigma(\triangle^{\frak n})$ of
$\triangle^{\frak n}$ is contained in
the spectrum $\sigma(\triangle^Z)$ of the differential form
Laplacian on $\Omega^*(Z)$.

Suppose that $\{\Gamma_i\}_{i=1}^\infty$ is a sequence of lattices in $N$
with quotients $Z_i = \Gamma_i \backslash N$ such that
$\lim_{i \rightarrow \infty} \diam(Z_i) = 0$. Then $\{Z_i\}_{i=1}^\infty$
obviously converges to a point, with bounded sectional
curvature in the collapse. 
We see that
there are eigenvalues of $\triangle^{Z_i}$ which
are constant in $i$, namely those which come from 
$\sigma(\triangle^{\frak n})$. 
By Proposition \ref{prop2} below, the other eigenvalues go
to infinity as $i \rightarrow \infty$. If $N$ is nonabelian then there
are positive eigenvalues of $\triangle^{Z_i}$ which are constant in $i$. 

In terms of Theorem \ref{th2}, 
$B$ is a point, 
$E^* = \Lambda^*({\frak n}^*)$ and $A^\prime = A^\prime_{[0]} =
d^{\frak n}$. This
shows that the term $A^\prime_{[0]}$ does appear in examples. In fact,
$A^\prime_{[0]} = 0$ if and only if $N$ is abelian.

By choosing different 
left-invariant metrics on $N$, we can make $\sigma(\triangle^{\frak n})$
arbitrarily close to zero while keeping the sectional
curvature bounded.  (In fact,
the sectional curvature goes to zero.) This is a special case of Theorem 
\ref{th7}.
We see that in general, 
there are no nontrivial lower bounds on the positive eigenvalues of
$\triangle^Z$ under the assumptions of bounded sectional curvature and 
bounded diameter.\\ \\
{\bf Example 2 :} Let $M$ be a compact 
manifold with a free $T^k$-action.
Let $g^{TM}$ be a $T^k$-invariant Riemannian metric on $M$. Then for
$\epsilon > 0$, there is a Riemannian metric $g^{TM}_\epsilon$ obtained
by multiplying $g^{TM}$ in the direction of the $T^k$-orbit by
$\epsilon$. Clearly $\lim_{\epsilon \rightarrow 0}(M, g^{TM}_\epsilon)
= M/T^k$, the collapse being with bounded sectional curvature 
\cite{Cheeger-Gromov (1986)}. This collapsing
is an example of the so-called adiabatic limit, for which the eigenvalues
of the differential form Laplacian have been studied in 
\cite{Berthomieu-Bismut (1994),Dai (1991),Forman (1995),Mazzeo-Melrose (1990)}.
Let $E$ be the flat ``cohomology'' vector bundle on $M/T^k$ with fiber
$\HH^*(T^k; \R)$; in fact, it is a trivial bundle. The results of
the cited references imply that as $\epsilon \rightarrow 0$, the
eigenvalues of $\triangle^M$ which remain finite
approach those of the Laplacian on
$\Omega^*(M/T^k; E)$. In particular,
the number of eigenvalues of the $p$-form Laplacian which go to zero
as $\epsilon \rightarrow 0$ is
$\sum_{a+b=p} \dim \left(
\HH^a(M/T^k; E^b)) \right)$, which is also the dimension of the
$E_2$-term of the Leray spectral sequence to compute
$\HH^p(M; \R)$.  This is consistent with Theorems \ref{th6} and \ref{th7}. 
Let $A^\prime_\epsilon$ be the superconnection on $E$ coming from Theorem
\ref{th2}, using $g^{TM}_\epsilon$. Then $\lim_{\epsilon \rightarrow 0}
A^\prime_\epsilon = \nabla^E$.\\ \\
{\bf Example 3 :} Suppose that $M$ is the total space of an oriented
circle bundle, with an $S^1$-invariant Riemannian metric. For $k \in \Z^+$,
consider the subgroup $\Z_k \subset S^1$. 
Then $\lim_{k \rightarrow \infty}
M/\Z_k = M/S^1$, the collapse obviously being with bounded sectional curvature.
By Fourier analysis, one finds that as $k \rightarrow \infty$, the
spectrum of $\triangle^{M/\Z_k}$ approaches the spectrum of the Laplacian
on $S^1$-invariant (not-necessarily-basic) differential forms on $M$.
In terms of Theorem \ref{th2}, $B = M/S^1$ and 
$E$ is the direct sum of two trivial $\R$-bundles on $B$. Let
$T$ be the curvature $2$-form of the fiber bundle 
$M \rightarrow M/S^1$. Then one finds that
the Laplacian acting on $S^1$-invariant forms
on $M$ is isomorphic to the Laplacian $\triangle^E = A^\prime \left(
A^\prime \right)^* + \left( A^\prime \right)^* A^\prime$, where
$A^\prime$ is the extension of the superconnection on   
$C^\infty (B;  E) = C^\infty(B) \oplus C^\infty(B)$ given by
\begin{equation} \label{eq2.1}
A^\prime \: = \: 
\begin{pmatrix}
\nabla^{E^0} & T \\
0 & \nabla^{E^1}
\end{pmatrix}.
\end{equation}
Here $\nabla^{E^0}$ and $\nabla^{E^1}$ are product connections. 
This shows that the term $A^\prime_{[2]}$ does appear in examples.
Note that if $M$ is simply-connected then 
$\{M/\Z_k\}_{k=1}^\infty$ are mutually nondiffeomorphic.

\section{Infranilmanifolds} \label{sect3}

In this section we first recall some basic facts
about infranilmanifolds.  Then in Proposition \ref{prop1} we show that the
orthogonal projection onto the parallel differential forms of $Z$ comes
from an averaging technique and so is independent of the choice of
parallel metric on $Z$, a result that will be crucial in what follows.

Let $N$ be a simply-connected connected nilpotent Lie group. 
Following \cite{Cheeger-Fukaya-Gromov (1992)}, when $N$ acts on a manifold
on the left we will denote it by $N_L$ and when it acts on a manifold
on the right we will denote it by $N_R$. As in
\cite{Cheeger-Fukaya-Gromov (1992)}, let us recall the elementary but
confusing point
that the right action of $N$ on $N$ generates left-invariant vector fields,
while the left action of $N$ on $N$ generates right-invariant vector fields.

There is 
a flat linear connection $\nabla^{aff}$ on $N$
which is characterized by the fact that
left-invariant vector fields are parallel. The group $\Aff(N)$ of
diffeomorphisms of $N$ which preserve $\nabla^{aff}$ is
isomorphic to $N_L \: \widetilde{\times} \: \Aut(N)$.

Suppose that $\Gamma$ is a discrete subgroup of $\Aff(N)$ which acts
freely and cocompactly on $N$, with $\Gamma \cap N_L$ of finite index in
$\Gamma$. Then the quotient space $Z = \Gamma \backslash
N$ is an infranilmanifold modeled on $N$.
We have the short exact sequences
\begin{equation} \label{eq3.1}
1 \longrightarrow N_L \longrightarrow \Aff(N)
\stackrel{p}{\longrightarrow} \Aut(N) \longrightarrow 1 
\end{equation}
and
\begin{equation} \label{eq3.2}
1 \longrightarrow \Gamma \cap N_L \longrightarrow \Gamma
\stackrel{p}{\longrightarrow} p(\Gamma) \longrightarrow 1. 
\end{equation}
Put $\widehat{\Gamma} = \Gamma \cap N_L$ and $F = p(\Gamma)$. Then $F$ is a
finite group. There is a normal cover 
$\widehat{Z} = \widehat{\Gamma} \backslash N$ of $Z$
with covering group $F$.

The connection $\nabla^{aff}$ descends to a flat 
connection on $TZ$, which we again
denote by $\nabla^{aff}$. Let $\Aff(Z)$ denote the affine group of $Z$, let
$\Aff_0(Z)$ denote the connected component of the identity in $\Aff(Z)$ and
let $\aff(Z)$ denote the affine Lie algebra of $Z$. Any element of
$\Aff(Z)$ can be lifted to an element of $\Aff(N)$. That is,
$\Aff(Z) = \Gamma \backslash (N_\Gamma {\Aff(N)})$, where 
$N_\Gamma {\Aff(N)}$ is the
normalizer of $\Gamma$ in $\Aff(N)$. Similarly,
$\Aff_0(Z) = C(\Gamma) \backslash (C_\Gamma {\Aff(N)})$, 
where $C_\Gamma {\Aff(N)}$ is the
centralizer of $\Gamma$ in $\Aff(N)$ and $C(\Gamma)$ is the center of
$\Gamma$.  
There is a short exact sequence
\begin{equation} \label{eq3.3}
1 \longrightarrow \Aff_0(Z) \longrightarrow 
\Aff({Z}) \longrightarrow \Out({\Gamma}) \longrightarrow 1.
\end{equation} 

As affine vector fields on $Z$ can be lifted to $F$-invariant affine
vector fields on $\widehat{Z}$, we have $\aff(Z) = \aff(\widehat{Z})^F$.
If $C(N)$ denotes the center of $N$ then
$\Aff_0(\widehat{Z}) = (\widehat{\Gamma} \cap C(N_R)) 
\backslash N_R$. In particular, if ${\frak n}$ is the Lie algebra of $N$
then $F$ acts by automorphisms on ${\frak n}$ and $\aff(Z) = {\frak n}_R^F$.

The $F$-invariant subspace $\Lambda^*({\frak n}^*)^F$ of
$\Lambda^*({\frak n}^*)$ is isomorphic to the vector space of 
differential forms on $Z$ which are parallel with
respect to $\nabla^{aff}$, or equivalently, to the 
$(N_L \widetilde{\times} F)$-invariant subspace of $\Omega^*(N)$.

Let $g^{TZ}$ be a Riemannian metric on $Z$ which is parallel with respect
to $\nabla^{aff}$. Such metrics correspond to $F$-invariant inner 
products on ${\frak n}$. Let $\diam(Z)$ denote the diameter of $Z$,
let $\nabla^Z$ denote the Levi-Civita connection of $Z$
and let $R^Z$ denote the Riemann curvature tensor of $Z$.

Let ${\cal P} : \Omega^*(Z) \rightarrow 
\Lambda^*({\frak n}^*)^F$ be orthogonal projection onto parallel
differential forms.

\begin{proposition} \label{prop1}
The orthogonal projection ${\cal P}$ 
is independent of the parallel metric $g^{TZ}$.
\end{proposition}
\begin{pf}
We first consider the case when $F = \{e\}$, so that
$Z$ is a nilmanifold $\Gamma \backslash N$.
As $N$ is nilpotent, it has a bi-invariant Haar measure $\mu$. We
normalize $\mu$ so that $\int_{\Gamma \backslash N} d\mu = 1$. 
Given
$\omega \in \Omega^*(Z)$, let
$\widetilde{\omega} \in \Omega^*(N)$ be its
pullback to $N$.  If $L_g$ denotes the left action of 
$g \in N_L$ on $N$ then for all $\gamma \in \Gamma$,
\begin{equation} \label{eq3.13}
L_{\gamma g}^* \: \widetilde{\omega} \: = \: 
L_{g}^* \: L_{\gamma}^* \:  \widetilde{\omega} \: = \: 
L_{g}^* \: \widetilde{\omega}.
\end{equation}
Hence it makes sense to define
$\overline{\widetilde{\omega}} \in \Omega^*(N)$ by
\begin{equation} \label{eq3.14}
\overline{\widetilde{\omega}} = 
\int_{{\Gamma} \backslash N_L}  \left( L_g^* \: \widetilde{\omega} 
\right) \: d\mu(g).
\end{equation}
For $h \in N_L$,
\begin{align} \label{eq3.15}
L_h^* \: \overline{\widetilde{\omega}} \:  = \: 
& \int_{{\Gamma} \backslash N_L}  \left( L_h^* \: L_g^* \: 
\widetilde{\omega} 
\right) \: d\mu(g) \: = \:
\int_{{\Gamma} \backslash N_L}  \left( L_{gh}^* \:
\widetilde{\omega} 
\right) \: d\mu(g) \\
= \: & \int_{{\Gamma} \backslash N_L}  \left( L_{g}^* \:
\widetilde{\omega} 
\right) \: d\mu(gh^{-1}) \: = \:
\int_{{\Gamma} \backslash N_L}  \left( L_{g}^* \:
\widetilde{\omega} 
\right) \: d\mu(g) \: = \: \overline{\widetilde{\omega}}. \notag
\end{align}
Thus $\overline{\widetilde{\omega}}$ is $N_L$-invariant and, in particular,
descends to a form $\overline{\omega} \in \Omega^*(Z)$.
Put
$P(\omega) = \overline{\omega}$. Then $P$ is idempotent, with $\Image(P)$
being the parallel differential forms. By construction,
$P$ is independent of the choice
of $g^{TZ}$. 
It remains to show that $P$ is self-adjoint. Given $\eta \in
\Omega^*(Z)$, let $\widetilde{\eta}$ be its lift to $N$. Consider the
function $f : N \times N \rightarrow \R$ given by 
\begin{equation} \label{eq3.16}
f(g, n) \: = \: \langle \widetilde{\eta}, L_g^* \: \widetilde{\omega} 
\rangle_n \: = \: \langle \widetilde{\eta}(n),  \widetilde{\omega}(gn) 
\rangle_n.
\end{equation}
For $\gamma \in \Gamma$, we have 
$f(\gamma g, n) = f(g \gamma^{-1}, \gamma n) = f(g, n)$.
It follows that we can write
\begin{equation} \label{eq3.17}
\langle \eta, P \omega \rangle_Z \: = \: \int_{(\Gamma \times \Gamma)
\backslash (N \times N)} \: 
\langle \widetilde{\eta}, L_g^* \: \widetilde{\omega} 
\rangle_n
 \: d\mu(g) \: d\mu(n),
\end{equation}
where the action of $\Gamma \times \Gamma$ on $N \times N$ is
$(\gamma_1, \gamma_2) \cdot (g, n) = (\gamma_1 g \gamma_2^{-1}, 
\gamma_2 n)$. Changing variable to $g^\prime \: = \: gn$, we have
\begin{align} \label{eq3.18}
\langle \eta, P\omega \rangle_Z \: & = \:   \int_{\Gamma 
\backslash N_L} \int_{\Gamma 
\backslash N_L} \: 
\langle \widetilde{\eta}, L_{g^\prime n^{-1}}^* \: \widetilde{\omega} 
\rangle_n
\: d\mu(g^\prime n^{-1}) \: d\mu(n) \\
& = \: 
\int_{\Gamma 
\backslash N_L} \int_{\Gamma 
\backslash N_L} \: 
\langle \widetilde{\eta}, L_{n^{-1}}^* L_{g^\prime}^* \: \widetilde{\omega} 
\rangle_n
\: d\mu(g^\prime) \: d\mu(n) \notag \\
& = \: \int_{\Gamma 
\backslash N_L} \int_{\Gamma 
\backslash N_L} \: 
\langle L_n^* \widetilde{\eta}, L_{g^\prime}^* \: \widetilde{\omega} 
\rangle_e
\: d\mu(g^\prime) \: d\mu(n) \notag \\
& = \: \int_{\Gamma 
\backslash N_L} \int_{\Gamma 
\backslash N_L} \: 
\langle L_{g^\prime}^* \: \widetilde{\omega}, L_n^* \widetilde{\eta}  
\rangle_e
\: d\mu(n) \: d\mu(g^\prime) \notag \\
& = \:
 \langle \omega, P\eta \rangle_Z
\: = \: \langle P\eta, \omega \rangle_Z. \notag
\end{align}
Thus $P$ is self-adjoint.

In the case of general $F$, we can apply the above argument equivariantly on
$\widehat{Z}$ with respect to $F$. 
As $F$ acts isometrically on $\widehat{Z}$, it commutes
with the orthogonal projection ${\cal P}$ on $\widehat{Z}$. 
As $F$ preserves $\mu$, it also commutes with the averaging operator
$P$ on $\widehat{Z}$. The proposition follows.
\end{pf}

\section{Eigenvalue Estimates on Infranilmanifolds} \label{infranil}

In this section we show, in Proposition \ref{prop2}, that if an
infranilmanifold $Z$ has
bounded sectional curvature and a diameter which goes to zero,
then all of the eigenvalues of $\triangle^Z$
go to infinity except for those that correspond
to eigenforms which are parallel on $Z$.

Let $N$ be a simply-connected connected $n$-dimensional
nilpotent Lie group with a
left-invariant Riemannian metric.
Let $\{e_i\}_{i=1}^{n}$  be an orthonormal basis
of ${\frak n}$.
Define the structure constants of
${\frak n}$ by $[e_i, e_j] \: = \: \sum_{k=1}^n c^k_{\: ij} \: e_k$.
Take the corresponding left-invariant basis $\{e_i\}_{i=1}^{n}$ of
$TN$, with dual basis of $1$-forms
$\{\tau^i \}_{i=1}^{n}$. Then the
components $\omega^i_{\: j} \: = \: \sum_{k=1}^n 
\: \omega^i_{\: jk} \: \tau^k$
of the Levi-Civita connection $1$-form 
$\omega \: = \: \sum_k \: \omega_k \: \tau^k$ are the constant matrices
\begin{equation} \label{eq3.22}
\omega^i_{\: jk} = - \: \frac{1}{2} \: \left( c^i_{\: jk} \: - \:
c^j_{\: ik} \: - \: c^k_{\: ij} \right).
\end{equation}

\begin{lemma} \label{lemma2}
Let $\kappa$ denote the scalar curvature of $Z$. Then
\begin{equation} \label{eq3.27}
\sum_{i,j,k = 1}^n \left( c^i_{\: jk} \right)^2 \: = \: - \: 4 \:
\kappa.
\end{equation}
\end{lemma}
\begin{pf}
As $\sum_{i,j,k = 1}^n \left( c^i_{\: jk} \right)^2$ is independent of the
choice of orthonormal basis, we will compute it using a special
orthonormal basis.
Recall the definition of
${\frak n}_{[k]}$ from (\ref{eq1.9}).
In particular, ${\frak n}_{[S]} \: = \: {\frak c}({\frak n})$. Following
the notation of \cite[\S 6]{Fukaya (1989)}
we take an orthonormal
basis $\{e_i\}_{i=1}^{n}$ of ${\frak n}$
such that $e_i \in {\frak n}_{[O(i)]}$ for some nondecreasing function
\begin{equation} \label{eq3.21}
O : \{1, \ldots, n\} \rightarrow \{0, \ldots S\},
\end{equation}
and $e_i \perp {\frak n}_{[O(i)+1]}$.

For a general Riemannian manifold,
we have the structure equations
\begin{align} \label{eq3.28}
d \tau^i \: & = \: - \: \sum_j \: \omega^i_{\: j} \: \wedge \: \tau^j, \\ 
\Omega^i_{\: j} \: & = \: d \omega^i_{\: j} \: + \: \sum_m \: 
\omega^i_{\: m} \: \wedge \: \omega^m_{\: \: j}. \notag
\end{align}
Then
\begin{align} \label{eq3.29}
\Omega^i_{\: j} \: & = \: d \: \sum_l \: \omega^i_{\: jl} \: \tau^l \: + \:
\sum_m \: \omega^i_{\: m} \: \wedge \: \omega^m_{\: \: j} \\
& = \: \sum_{k,l} \: \left( e_k \omega^i_{\: jl} \right) 
\: \tau^k \: \wedge \tau^l \: + 
\: \sum_m \: \omega^i_{\: jm} \: d \tau^m \: +
\sum_{k,l,m} \:
\omega^i_{\: mk} \: \omega^m_{\: \: jl} \: \tau^k \: \wedge \tau^l.  \notag
\end{align}
This gives the Riemann curvature tensor as
\begin{equation} \label{eq3.30}
R^i_{\: jkl} \: = \: e_k \omega^i_{\: jl} \: - \: e_l \omega^i_{\: jk} \:
+ \: \sum_m \left[- \: \omega^i_{\: jm} \: \omega^m_{\: \: lk}
+ \: \omega^i_{\: jm} \: \omega^m_{\: \: kl}
+ \: \omega^i_{\: mk} \: \omega^m_{\: \: jl}
- \: \omega^i_{\: ml} \: \omega^m_{\: \: jk} \right].
\end{equation}
Then
\begin{align} \label{eq3.31}
\kappa \: & = \: \sum_{i,j} \left(  e_i \omega^i_{\: jj} \: - \: 
e_j \omega^i_{\: ji}  \right) \:
+ \: \sum_{i,j,m} \left[- \: \omega^i_{\: jm} \: \omega^m_{\: \: ji}
+ \: \omega^i_{\: jm} \: \omega^m_{\: \: ij}
+ \: \omega^i_{\: mi} \: \omega^m_{\: \: jj}
- \: \omega^i_{\: mj} \: \omega^m_{\: \: ji} \right] \notag \\
& = \: \sum_{i,j} \left(  e_i \omega^i_{\: jj} \: - \: 
e_j \omega^i_{\: ji}  \right) \:
+ \: \sum_{i,j,m} \left[\omega^i_{\: jm} \: \omega^m_{\: \: ij}
+ \: \omega^i_{\: mi} \: \omega^m_{\: \: jj} \right].
\end{align}

In our case, the components of the connection matrix are constant.
Also, as ${\frak n}$ is nilpotent,
\begin{equation} \label{eq3.32}
\omega^i_{\: jj} \: = \: c^j_{\: ij} \: = 0.
\end{equation}
Then one obtains
\begin{equation} \label{eq3.33}
\kappa \: = \: - \:
 \sum_{i,j,k = 1}^n 
\: \omega^i_{\: jk} \: \omega^i_{\: kj}.
\end{equation}
Separating $\omega^i_{\: jk}$ into its components which are symmetric or
antisymmetric in $j$ and $k$, and using (\ref{eq3.22}), we obtain
\begin{align} \label{eq3.34}
\kappa \: & = \: - \: \frac{1}{4} \:
 \sum_{i,j,k=1}^n 
\: \left( c^j_{\: ik} \: + \: c^k_{\: ij}\right)^2 \: + \: \frac{1}{4}
\:  \sum_{i,j,k=1}^n \: \left( c^i_{\: jk} \right)^2 \\
& = \: - \: \sum_{i,j,k=1}^n  \: \left[ \frac{1}{2} \:
\: c^j_{\: ik} \: c^k_{\: ij} \: + \: \frac{1}{4}
\: \left( c^i_{\: jk} \right)^2 \right] \notag.
\end{align}
As ${\frak n}$ is nilpotent, it follows that $c^j_{\: ik} \: c^k_{\: ij} \:
= \: 0$. This proves the lemma.
\end{pf}

Let $Z$ be an infranilmanifold with an affine-parallel metric.
Let $\triangle^Z$ denote
the Laplacian acting on $\Omega^*(Z)$. Let $\triangle^{inv}$ be the
finite-dimensional Laplacian acting on $\Lambda^*({\frak n}^*)^F$. 
\begin{proposition} \label{prop2}
There are positive constants $A$ and $A^\prime$, depending only on
$\dim(Z)$, such that if $\parallel R^Z \parallel_\infty  \diam(Z)^2 \: \le 
A^\prime$ then
the spectrum $\sigma(\triangle^Z)$ of $\triangle^Z$ satisfies
\begin{equation} \label{eq3.19}
\sigma(\triangle^Z) \cap \left[ 0, \: A \: 
{\diam(Z)^{-2}} \right) =
\sigma(\triangle^{inv}) \cap \left[ 0, \: A \:
{\diam(Z)^{-2}} \right).
\end{equation}
\end{proposition}
\begin{pf}
Recall the definition of ${\cal P}$ from Proposition \ref{prop1}.
It is enough to show that under the hypotheses of the present proposition,
the spectrum of $\triangle^Z$ on
$\Ker({\cal P})$ is bounded below by
$A \: \diam(Z)^{-2}$.

As $\widehat{Z}$ isometrically covers $Z$ with covering group $F$,
the spectrum of $\triangle^Z$ on $\Ker({\cal P}) \subset
\Omega^*(Z)$ is contained in the spectrum of 
$\triangle^{\widehat{Z}}$ on $\Ker({\cal P}) \subset
\Omega^*(\widehat{Z})$. 

\begin{lemma}
There is a function $\eta : \N \rightarrow \N$ such that
\begin{equation} \label{eq3.20}
\diam(\widehat{Z}) \: \le \: \eta (|F|) \: \diam({Z}).
\end{equation}
\end{lemma}
\begin{pf}
Let $\widehat{z}_1, \widehat{z}_2 \in \widehat{Z}$ be such that
$\diam(\widehat{Z}) = d(\widehat{z}_1, \widehat{z}_2)$. It is easy to
see that $d(\widehat{z}_1, F \cdot \widehat{z}_2) \le \diam(Z)$.
Let $z_2 \in Z$ be the projection of $\widehat{z}_2 \in \widehat{Z}$.
Then it is enough to bound $\parallel \cdot \parallel_{geo}$ from above 
on $\pi_1(Z, z_2) \cong F$,
i.e. to bound the minimal lengths of curves in the classes of
$\pi_1(Z, z_2)$. From
\cite[Proposition 3.22]{Gromov (1999)}, there is a set of generators
of $\pi_1(Z, z_2)$ on which $\parallel \cdot \parallel_{geo}$ is bounded above
by $2 \: \diam(Z)$. Given
$r \in \N$, 
there is a finite number of groups of order $r$, up to isomorphism,
and each of these groups has a finite number of generating sets. The lemma 
follows.
\end{pf}

Furthermore, there is a universal bound
$|F| \le \const (\dim(Z))$ \cite{Buser-Karcher (1981)}.
Hence without loss of generality, we may assume that
$F = \{e\}$ so that $Z$ is a nilmanifold $\Gamma \backslash N$. 

Let $E^i$ denote exterior multiplication on $\Omega^*(Z)$ by
$\tau^i$ and let $I^i$ denote interior multiplication by
$e_i$. 
From the Bochner formula, if $\eta \in \Omega^*(Z)$ then
\begin{equation} \label{add1}
\langle \eta, \triangle^{Z} \eta \rangle_{Z} \: = \:
\langle \nabla^Z \eta, \nabla^Z \eta \rangle_Z \: + \:
\sum_{ijkl} \:
\int_{Z} R^{Z}_{ijkl} \langle E^i \: I^j \: \eta, E^k \: I^l \: 
\eta \rangle
\: d\vol_{Z}.
\end{equation}
Using the left-invariant vector fields on $N$, there is an 
isometric isomorphism
\begin{equation} \label{add2}
\Omega^*(Z) \cong C^\infty(Z) \: \otimes \: 
\Lambda^*({\frak n}^*).
\end{equation}
With respect to this isomorphism,
\begin{equation} \label{add3}
\nabla^{TZ}_{e_i} \: = \: \left( e_i \otimes \Id \right) \: + \:
\left( \Id \otimes 
\sum_{j,k} \: \omega^j_{\: ki} \: E^j \: I^k \right),
\end{equation}
where $E^j$ and $I^k$ now act on $\Lambda^*({\frak n}^*)$.
It follows that
\begin{equation} \label{add3.25}
\langle \nabla^Z \eta, \nabla^Z \eta \rangle_Z \: \ge \:
\sum_i \langle (e_i \otimes \Id) \eta, (e_i \otimes \Id)\eta \rangle_Z
\: - \: \sum_i \big|
\sum_{j,k} \: \omega^j_{\: ki} \: E^j \: I^k \eta
\big|_Z^2.
\end{equation}
Let $\triangle^Z_0$ be the ordinary Laplacian on $C^\infty(Z)$. 
With respect to (\ref{add2}), consider the operator
$\triangle^Z_0 \: \otimes \: \Id$. We have
\begin{equation} \label{add3.5}
\langle \eta, (\triangle^{Z}_0 \: \otimes \: \Id
) \eta \rangle_{Z} \: = \: 
\sum_i \langle (e_i \otimes \Id) \eta, (e_i \otimes \Id)\eta \rangle_Z
\end{equation}
Using (\ref{eq3.22}),
(\ref{add1}), (\ref{add3.25}), (\ref{add3.5}) and 
Lemma \ref{lemma2}, we obtain
\begin{equation} \label{add4}
\langle \eta, \triangle^{Z} \eta \rangle_{Z} \: \ge \:
\langle \eta, (\triangle^{Z}_0 \: \otimes \: \Id
) \eta \rangle_{Z} \: - \: \const
\parallel R^Z \parallel_\infty \: | \eta|_Z^2.
\end{equation}

In terms of (\ref{add2}), $\Ker({\cal P}) \: \cong \: 1^\perp \: \otimes \:
\Lambda^*({\frak n}^*)$, where $1$ denotes the constant function on $Z$.
Thus if $\eta \in \Ker({\cal P})$ then
$\langle \eta, (\triangle^{Z}_0 \: \otimes \: \Id)
\eta \rangle_{Z} \: \ge \: \lambda_{0,2}
\: | \eta|_Z^2$, where $\lambda_{0,2}$ is the first positive eigenvalue
of the function Laplacian on $Z$. There is a lower bound
\begin{equation}
\lambda_{0,2} \: \ge \: \diam(Z)^{-2} \: 
f \left( \parallel R^Z \parallel_\infty \: \diam(Z)^2 \right)
\end{equation}
for some
smooth function $f$ with $f(0) > 0$
\cite{Berard (1988)}. Thus the spectrum of $\triangle^Z$ on $\Ker({\cal P})$
is bounded below by
\begin{equation}
\diam(Z)^{-2} \: \left[
f \left( \parallel R^Z \parallel_\infty \: \diam(Z)^2 \right) \: - \:
\const \: \parallel R^Z \parallel_\infty \: \diam(Z)^2 \right].
\end{equation}
Taking $A \: = \: \frac{3}{4} \: f(0)$, the proposition follows. 
\end{pf}

\section{Affine Fiber Bundles} \label{sect4}

In this section we first show that the eigenvalues of a superconnection
Laplacian are continuous with respect to the superconnection, the Riemannian
metric and the Euclidean inner product.
We then construct the superconnection $A^\prime$ associated to
an affine fiber bundle $M \rightarrow B$ and prove Theorem \ref{th1}.

Let $B$ be a smooth connected closed Riemannian manifold.
Let $E = \oplus_{j=0}^m E^j$ be a 
$\Z$-graded real vector bundle on $B$. For background information about 
superconnections, we refer to 
\cite[Chapter 1.4]{Berline-Getzler-Vergne (1992)}, 
\cite{Bismut (1986)},
\cite{Bismut-Lott (1995)} and \cite{Quillen (1985)}.
Let $A^\prime$ be a degree-$1$ superconnection on $E$.  That is,
$A^\prime$ is an $\R$-linear map from $C^\infty(B; E)$ to $\Omega (B; E)$
with a decomposition
\begin{equation} \label{eq4.1}
A^\prime \: = \: \sum_{k=0}^{\dim(B)} A^\prime_{[k]}
\end{equation}
where
\begin{itemize}
\item $A^\prime_{[1]}$ 
is a connection $\nabla^E$ on $E$ which preserves the $\Z$-grading.
\item For $k \ne 1$, 
$A^{\prime}_{[k]} \in \Omega^k \left(B; \Hom(E^*, E^{*+1-k})
 \right)$.
\end{itemize}
We can extend $A^\prime$ to an $\R$-linear map on $\Omega (B; E)$ using the
Leibniz rule.  We assume that $A^\prime$ is flat, in the sense that
\begin{equation} \label{eq4.2}
(A^\prime)^2 = 0.
\end{equation}
Let $h^E$ be a Euclidean inner product on $E$ such that $E^j$ is orthogonal to
$E^{j^\prime}$ if $j \ne j^\prime$. Let $(A^\prime)^*$ be the adjoint
superconnection with respect to $h^E$ and put 
\begin{equation} \label{eq4.3}
\triangle^E \: = \: A^\prime (A^\prime)^* \: + \:
(A^\prime)^* \: A^\prime.
\end{equation}
Then $\triangle^E$ preserves the total $\Z$-grading on $\Omega (B; E)$ and
decomposes with respect to the grading as 
$\triangle^E = \bigoplus_{p} \triangle^E_p$. By elliptic theory,
$\triangle^E_p$ has a discrete spectrum.

If $g^{TB}_1$ and $g^{TB}_2$ 
are two Riemannian metrics on $B$ and $\epsilon \ge 0$, 
we say that $g^{TB}_1$ and $g^{TB}_2$ are $\epsilon$-close if
\begin{equation} \label{eq4.4}
e^{-\epsilon} \: g^{TB}_2 \: \le \: g^{TB}_1 \: \le \: e^\epsilon \: 
g^{TB}_2.
\end{equation}
Similarly, if $h^E_1$ and $h^E_2$ are two Euclidean inner products on $E$,
we say that $h^E_1$ and $h^E_2$ are $\epsilon$-close if
\begin{equation} \label{eq4.5}
e^{-\epsilon} \: h^E_2 \: \le \: h^E_1 \: \le \: e^\epsilon \: h^E_2.
\end{equation}
If $S_1 = \{ \lambda_{1,j} \}$ and 
$S_2 = \{ \lambda_{2,j} \}$ are two countable nondecreasing ordered sets of
nonnegative real
numbers then we say that $S_1$ and $S_2$ are $\epsilon$-close if for
all $j$,
\begin{equation} \label{eq4.6}
e^{-\epsilon} \: \lambda_{2,j} \: \le \: \lambda_{1,j} \: \le \: e^\epsilon \: 
\lambda_{2,j}.
\end{equation}

For simplicity, we will omit the subscript $p$, the form degree, in this
section when its role is obvious.

\begin{lemma} \label{lemma4}
There is an integer $J = J(\dim(B)) > 0$ such that if $g^{TB}_1$ and 
$g^{TB}_2$ are $\epsilon$-close,
and $h^E_1$ and $h^E_2$ are $\epsilon$-close, then the corresponding Laplacians
$\triangle^E_1$ and $\triangle^E_2$ have spectra which are $J \epsilon$-close.
\end{lemma}
\begin{pf}
As in \cite[Prop. 3.1]{Dodziuk (1982)}, using a trick apparently first
due to Cheeger, 
we can write the
spectrum of $\triangle^E$ on $\Image \left( (A^\prime)^* \right)$ as
\begin{equation} \label{eq4.7}
\lambda_j \: = \: \inf_{V} \sup_{\eta \in V - \{0\}} 
\sup_{\theta \in \Omega(B; E)}
\left\{ \frac{\langle \eta,\eta \rangle}{\langle \theta,\theta \rangle} : 
\eta = A^\prime \theta \right\},
\end{equation}
where $V$ ranges over $j$-dimensional subspaces of $\Image(A^\prime)$.
As the Riemannian metric and Euclidean inner
product only enter in defining $\langle \cdot, \cdot \rangle$, 
the lemma follows as in
\cite{Dodziuk (1982)}.
\end{pf}

We will also need a result about how the spectrum of $\triangle^E$ depends
on the superconnection $A^\prime$. 
Given
$X \in \Omega (B; \End(E))$, let $\parallel X \parallel$ be the 
operator norm for the action of $X$ on the $L^2$-completion of
$\Omega (B; E)$.
If $A_1^\prime$ and $A_2^\prime$ are two superconnections as above
then $A_1^\prime - A_2^\prime \in \Omega (B; \End(E))$.  
Fix $g^{TB}$ and $h^E$.

\begin{lemma} \label{lemma?}
For all $j \in \Z^+$,
\begin{equation} \label{eq5.13}
|\lambda_j(A^\prime_1)^{1/2} - \lambda_j(A^\prime_2)^{1/2}| \: \le \:
(2 + \sqrt{2}) \: \parallel A_1^\prime - A_2^\prime \parallel.
\end{equation}
\end{lemma}
\begin{pf}
Put $x \: = \: \parallel A_1^\prime - A_2^\prime \parallel$.
If $\omega \in \Omega (B; E)$ is nonzero then
\begin{equation} \label{eq5.14}
\left| \frac{|A_1^\prime \omega|}{|\omega|} - 
\frac{|A_2^\prime \omega|}{|\omega|} \right| \: \le 
\frac{|(A_1^\prime - A_2^\prime) \omega|}{|\omega|} \: \le 
\: x
\end{equation}
and
\begin{equation} \label{eq5.15}
\left| \frac{|(A_1^\prime)^* \omega|}{|\omega|} - 
\frac{|(A_2^\prime)^* \omega|}{|\omega|} \right| \: \le \:
\frac{|((A_1^\prime)^* - 
(A_2^\prime)^*) \omega|}{|\omega|} \: 
\le \: x
\end{equation}
Define $\vec{v}_1, \vec{v}_2 \in \R^2$ by
\begin{equation} \label{eq5.16}
\vec{v}_i \: = \: \left( \frac{|A_i^\prime \omega|}{|\omega|},
\frac{|(A_i^\prime)^* \omega|}{|\omega|} \right).
\end{equation}
Then (\ref{eq5.15}) and (\ref{eq5.16}) imply that 
\begin{equation} \label{eq5.17}
\left| \left(
\parallel \vec{v}_2 \parallel \: - \: \parallel \vec{v}_1 \parallel  
\right) \right| \: \le \: 
\parallel \vec{v}_2 - \vec{v}_1 \parallel \: \le \: \sqrt{2} \: x.
\end{equation}
Hence
\begin{equation} \label{eq5.18}
\parallel \vec{v}_2 \parallel^2 - \parallel \vec{v}_1 \parallel^2 
 \: = \: 
\left( \parallel \vec{v}_2 \parallel - \parallel \vec{v}_1 \parallel 
\right) \cdot \left(
\parallel \vec{v}_2 \parallel + \parallel \vec{v}_1 \parallel \right) 
\: \le \: \sqrt{2} \: x \: \left( 2 \parallel \vec{v}_2 \parallel \: + \:
\sqrt{2} \: x \right), 
\end{equation} 
so
\begin{align} \label{eq5.19}
\parallel \vec{v}_1 \parallel^2 \: &\ge \: \parallel \vec{v}_2 \parallel^2 
\: - \: 2 \: \sqrt{2} \: x \: \parallel \vec{v}_2 \parallel \: - \: 2 \:
 x^2 \\
& = \: \left(  \parallel \vec{v}_2 \parallel 
\: - \: \sqrt{2} \: x \right)^2 \: - \: 4 \: x^2. \notag
\end{align}
Thus
\begin{equation} \label{eq5.20}
\parallel \vec{v}_1 \parallel^2 \: \ge 
\: \max \left( 0, \left(  \parallel \vec{v}_2 \parallel 
\: - \: \sqrt{2} \: x \right)^2 \: - \: 4 \: x^2 \right),
\end{equation}
or equivalently,
\begin{equation} \label{eq5.21}
\frac{\langle \omega, \triangle_{A_1^\prime} \omega 
\rangle}{\langle \omega, \omega \rangle} \ge 
\: \max \left( 0, \left(  \left( \frac{\langle \omega, 
\triangle_{A_2^\prime} \omega 
\rangle}{\langle \omega, \omega \rangle} \right)^{1/2}
\: - \: \sqrt{2} \: x \right)^2 \: - \: 4 \: x^2 \right),
\end{equation}
The minmax characterization of
eigenvalues
\begin{equation} \label{eq5.22}
\lambda_j(A^\prime) \: = \: \inf_{V} \sup_{\omega \in V - \{0\}} 
\left\{ \frac{\langle \omega, \triangle_{A^\prime} \omega 
\rangle}{\langle \omega, \omega \rangle} \right\},
\end{equation}
where $V$ ranges over $j$-dimensional subspaces of $\Omega (B; E)$,
implies
\begin{equation} \label{eq5.23}
\lambda_j(A^\prime_1) \: \ge \: \max \left(0, 
\left( \lambda_j^{1/2}(A^\prime_2) -
\sqrt{2} \: x \right)^2 \: - \: 4 \: x^2 \right).
\end{equation}
An elementary calculation then gives
\begin{equation} \label{eq5.24}
\lambda_j(A^\prime_1)^{1/2} \: - \: \lambda_j(A^\prime_2)^{1/2} \: \ge \: - \:
(2 + \sqrt{2}) \: x. 
\end{equation}
Symmetrizing in $A^\prime_1$ and $A^\prime_2$, the proposition follows.
\end{pf}

Let $M$ be a closed manifold which is the total space of an affine
fiber bundle, as in Definition \ref{def2}. 
Let $T^HM$ be a horizontal
distribution on $M$ so that the corresponding holonomy on $B$ lies in 
$\Aff(Z)$. If $m \in Z_b$ then using $T^HM$, we can write 
$\Lambda^*(T_m^* M) \cong \Lambda^*(T_b^* B) \: \widehat{\otimes} \:
\Lambda^*(T_m^* Z_b)$. That is, we can compose differential forms on 
$M$ into their horizontal and vertical components.
Correspondingly, 
there is an infinite-dimensional $\Z$-graded real vector bundle $W$ on $B$
such that $\Omega^*(M) \cong \Omega (B; W)$; 
see \cite[Section
III(a)]{Bismut-Lott (1995)}. A fiber $W_b$ of $W$ is isomorphic
to $\Omega^*(Z_b)$. We will call $C^\infty(B; W)$ the vertical
differential forms.
The exterior derivative $d^M \: : \:
\Omega^*(M) \rightarrow \Omega^*(M)$, when considered to be an operator
$d^M \: : \:
\Omega (B; W) \rightarrow \Omega(B; W)$, is the extension to 
$\Omega (B; W)$ of
a flat degree-$1$ superconnection on $W$.
From \cite[Proposition 3.4]{Bismut-Lott (1995)}, we can write the 
superconnection as
\begin{equation} \label{eq4.8}
d^Z \: + \: \nabla^W \: + \: i_T,
\end{equation}
where
\begin{itemize}
\item  $d^Z \in C^\infty(B; \Hom(W^*; W^{*+1}))$ is vertical differentiation,
\item $\nabla^W \: : \: C^\infty(B; W) \rightarrow \Omega^1(B; W)$ comes from 
Lie differentiation in the horizontal direction and
\item 
$i_T \in \Omega^2(B; \Hom(W^*; W^{*-1}))$ is interior multiplication by the
curvature $2$-form $T \in \Omega^2(M; TZ)$ of $T^HM$.
\end{itemize}
Acting on $\Omega^*(M)$, we have
\begin{equation} \label{eq4.9}
d^M \: = \: d^Z \: + \: d^W \: + \: i_T,
\end{equation}
where $d^W \: : \: \Omega^*(B; W) \rightarrow \Omega^{*+1}(B; W)$ 
is exterior differentiation on $B$ using $\nabla^W$.

Let $E$ be the finite-dimensional subbundle
of $W$ such that $E_b$ consists of the elements of $\Omega^*(Z_b)$
which are parallel on $Z_b$.
The fibers of $E$ are isomorphic to $\Lambda^*({\frak n}^*)^F$
and $C^\infty(B; E)$ is isomorphic to
the vertical differential forms on $M$ whose restrictions to the fibers are
parallel. 
Furthermore, the superconnection (\ref{eq4.8}) restricts to
a flat degree-$1$ superconnection $A^\prime$ on $E$, as
exterior differentiation on $M$ preserves the space of
fibrewise-parallel differential forms. From (\ref{eq4.8}),
\begin{equation} \label{eq4.10}
A^\prime \: = \: d^{\frak n} \: + \: \nabla^E \: + \: i_T,
\end{equation}
where $d^{\frak n}$ is the differential on $\Lambda^*({\frak n}^*)^F$ and
$\nabla^E$ comes from $T^HM$ through the action of $\Aff(Z)$ on
$\Lambda^*({\frak n}^*)^F$. Acting on $\Omega (B; E)$, we have
\begin{equation} \label{eq4.11}
A^\prime \: = \: d^{\frak n} \: + \: d^E \: + \: i_T,
\end{equation}
where $d^E$ is exterior differentiation on $\Omega(B; E)$ 
using $\nabla^E$.\\ \\
{\bf Remark :} 
The connection $\nabla^E$ is generally not flat.
As $A^\prime$ is flat, we have
\begin{equation} \label{eq4.12}
\left( \nabla^E \right)^2 \: = \: - \left( 
d^{\frak n} \: i_T \: + \: i_T \: d^{\frak n} \right).
\end{equation}
Thus the curvature of $\nabla^E$ is given by Lie differentiation with respect
to the
(negative of the) curvature $2$-form $T$. More geometrically, given $b \in B$,
let $\gamma$ be a loop
in $B$ starting from $b$ and let $h(\gamma) \in \Aff(Z_{b})$ be the 
holonomy of the connection $T^HM$
around $\gamma$. Then the holonomy of $\nabla^E$ around
$\gamma$ is the action of $h(\gamma)$ on the fiber $E_{b}$. In particular,
the infinitesimal holonomy of $\nabla^E$ lies in the
image of the Lie algebra $\aff(Z)$ in $\End(E_{b})$. From the discussion
after (\ref{eq3.3}),
$\aff(Z)$ lies in ${\frak n}_R$. As the elements of $E_{b}$ are
$N_L$-invariant forms on $N$, they are generally not annihilated by
$\aff(Z)$. \\

Suppose in addition
that $M$ is a Riemannian affine fiber bundle, as in Definition 
\ref{def2}.
Then $g^{TZ}$ induces an $L^2$-inner product 
$h^W$ on $W$ and a Euclidean inner 
product
$h^E$ on $E$. Let $\diam(Z)$ denote the maximum diameter of the fibers 
$\{Z_b\}_{b \in B}$ in the intrinsic metric and let $\Pi$ denote the
second fundamental forms of the fibers.
From $g^{TZ}$, $T^HM$ and $g^{TB}$, we obtain a
Riemannian metric $g^{TM}$ on $M$. 
Let $\triangle^M$ denote
the Laplacian acting on $\Omega^*(M)$ and define
$\triangle^E$, acting on $\Omega (B; E)$, as in (\ref{eq4.3}).
Let $R^M$ denote the Riemann curvature
tensor of $g^{TM}$.

Let ${\cal P}^{fib}$ 
be fiberwise orthogonal projection from $\Omega (B; W)$ to 
$\Omega (B; E)$. We claim that 
${\cal P}^{fib}$ commutes with $d^M$. Looking at
(\ref{eq3.14}), ${\cal P}^{fib}$
clearly commutes with $d^Z$. Using the fact that the holonomy of $T^HM$ lies
in $\Aff(Z)$, it follows from (\ref{eq3.14}) and the proof of Proposition 
\ref{prop1} that ${\cal P}^{fib}$ commutes with $\nabla^W$.
As $T$ takes values in parallel vector fields on $Z$,
it follows from (\ref{eq3.14}) and the proof of Proposition 
\ref{prop1} that ${\cal P}^{fib}$ commutes with
$i_T$. Thus ${\cal P}^{fib}$
commutes with $d^M$. As the fiberwise metrics are parallel on the fibers, 
it follows that ${\cal P}^{fib}$ also commutes with $\left(d^M\right)^*$.

Then with respect to
the decomposition $\Omega^*(M) \: = \: \Image({\cal P}^{fib}) 
\: \oplus \: \Ker({\cal P}^{fib})$,
$\triangle^M$ is isomorphic to $\triangle^E \: \oplus \:
\triangle^M \big|_{Ker({\cal P}^{fib})}$. \\ \\
{\bf Proof of Theorem \ref{th1} } :
From Proposition \ref{prop2}, there is a constant $A \: > \: 0$ so that
for all $b \in B$, the spectrum of $\triangle^{Z_b} \big|_{Ker({\cal P})}$
is bounded below by $A \cdot \diam(Z_b)^{-2}$.
It suffices to show that there is a constant $C$ as in the statement of
the theorem such that
\begin{equation}
\sigma \left(\triangle^M \big|_{Ker({\cal P}^{fib})} \right) 
\subset \left[A \: 
{\diam(Z)^{-2}} \: - C \: \left(
\parallel R^{M} \parallel_\infty + 
\parallel \Pi \parallel_\infty^2 +
\parallel T \parallel_\infty^2 
 \right), \infty \right).
\end{equation}

We use the notation of \cite[Section III(c)]{Bismut-Lott (1995)} to
describe the geometry of the fiber bundle $M$.  In particular,
lower case Greek indices refer to horizontal directions, lower case
italic indices refer to vertical directions and upper case italic indices
refer to either.
Let $\{\tau^i\}_{i=1}^{dim(Z)}$ and 
$\{\tau^\alpha\}_{\alpha=1}^{dim(B)}$ be a local orthonormal basis of $1$-forms
as in \cite[Section III(c)]{Bismut-Lott (1995)}, with dual basis
$\{e_i\}_{i=1}^{dim(Z)}$ and 
$\{e_\alpha\}_{\alpha=1}^{dim(B)}$. Let $E^J$ be exterior multiplication
by $\tau^J$ and let $I^J$ be interior multiplication by $e_J$.
The tensors $\Pi$ and $T$ are parts of the connection
$1$-form component $\omega^i_{\: \alpha} \: = \: \sum_j \:
\omega^i_{\: \alpha j} \:
\tau^j \: + \: \sum_\beta \: \omega^i_{\: \alpha \beta} \: \tau^\beta$,
with symmetries
\begin{equation} \label{eq4.13}
\omega_{\alpha k j} \: = \: 
\omega_{\alpha j k} \: = \: 
- \: \omega_{j \alpha k} \: = \: 
- \: \omega_{k \alpha j},
\end{equation}
\begin{equation*}
\omega_{\beta \alpha j} \: = \: 
- \: \omega_{\alpha \beta j} \: = \: 
- \: \omega_{\alpha j \beta} \: = \: 
\omega_{j \alpha \beta} \: = \: 
\omega_{\beta j \alpha} \: = \: 
- \: \omega_{j \beta \alpha}. 
\end{equation*}

Given $\eta \in \Omega^*(M)$, the Bochner formula gives
\begin{equation} \label{eq4.14}
\langle \eta, \triangle^M \eta \rangle_M \: = \:
\langle \nabla^M \eta, \nabla^M \eta \rangle_M \: + \: \sum_{PQRS} \:
\int_M R^M_{PQRS} \langle E^P \: I^Q \: \eta, E^R \: I^S \: \eta \rangle
\: d\vol_M.
\end{equation}
Here 
\begin{equation} \label{eq4.15}
\nabla^{M} : C^\infty(M; \Lambda^* T^*M) \rightarrow
C^\infty(M; T^*M \otimes \Lambda^* T^*M)
\end{equation}
is, of course, the Levi-Civita connection on $M$. 
We can write $\nabla^M = \nabla^V + \nabla^H$ where 
\begin{equation} \label{eq4.16}
\nabla^{V} : C^\infty(M; \Lambda^* T^*M) \rightarrow
C^\infty(M; T^{vert}M \otimes \Lambda^* T^*M)
\end{equation}
denotes
covariant differentiation in the vertical direction and 
\begin{equation} \label{eq4.17}
\nabla^{H} : C^\infty(M; \Lambda^* T^*M) \rightarrow
C^\infty(M; T^{hor}M \otimes \Lambda^* T^*M)
\end{equation} denotes
covariant differentiation in the horizontal direction.  Then
\begin{align} \label{eq4.18}
\langle \eta, \triangle^M \eta \rangle_M \: & = \:
\langle \nabla^V \eta, \nabla^V \eta \rangle_M \: + \:
\langle \nabla^H \eta, \nabla^H \eta \rangle_M \: + \:
\int_M \sum_{PQRS} \:
R^M_{PQRS} \langle E^P \: I^Q \: \eta, E^R \: I^S \: \eta \rangle
\: d\vol_M. \\
& \ge \:
\langle \nabla^V \eta, \nabla^V \eta \rangle_M \: - \: \const \:
\parallel R^M \parallel_\infty \: \langle \eta, \eta \rangle_M \notag \\
& = \: \int_B \int_{Z_b} \left[
\left| \nabla^V \eta \right|^2(z) - \: \const \:
\parallel R^M \parallel_\infty \: |\eta(z)|^2 \right]
 \: d\vol_{Z_b} \:
d\vol_B. \notag
\end{align}

Let 
\begin{equation} \label{eq4.19}
\nabla^{TZ} : C^\infty(M; \Lambda^*(T^{*,vert}M)) \rightarrow 
C^\infty(M; T^*M \otimes \Lambda^*(T^{*,vert}M))
\end{equation}
denote the Bismut connection acting on  $\Lambda^*(T^{*,vert}M)$
\cite[Proposition 10.2]{Berline-Getzler-Vergne (1992)},
\cite[Definition 1.6]{Bismut (1986)}. 
On a given fiber $Z_b$, there is a canonical flat connection on $T^{hor}M
\big|_{Z_b}$. Hence we can use $\nabla^{TZ}$ to
vertically differentiate sections of $\Lambda^*(T^*M) =
\Lambda^*(T^{*,vert}M) \: \widehat{\otimes} \: \Lambda^*(T^{*,hor}M)$.
That is, we can define
\begin{equation} \label{eq4.20}
\nabla^{TZ} : C^\infty(M; \Lambda^*(T^*M)) \rightarrow 
C^\infty(M; T^{*,vert}M \otimes \Lambda^*(T^*M))
\end{equation}
Explicitly, with respect to a local framing,
\begin{equation} \label{eq4.21}
\nabla^{TZ}_{e_i} \eta \: = \: e_i \eta \: + 
\sum_{j,k} \: \omega^j_{\: ki} \: E^j \: I^k \: \eta
\end{equation}
and
\begin{equation} \label{eq4.22}
\nabla^V_{e_i} \eta \: = \nabla^{TZ}_{e_i} \eta
\: + \: \sum_{j \alpha} \: \omega^j_{\: \alpha i} \: E^j \: I^\alpha \:
\eta \: + \: \sum_{\alpha k} \: \omega^\alpha_{\: ki} \: E^\alpha \: I^k \:
\eta \: + \: \sum_{\alpha \beta} 
\: \omega^\alpha_{\: \beta i} \: E^\alpha \: I^\beta \: \eta.
\end{equation} 
Then from (\ref{eq4.21}) and (\ref{eq4.22}),
\begin{equation} \label{eq4.23}
\langle \nabla^V \eta, \nabla^V \eta \rangle_M \: \ge \:
\int_B \int_{Z_b} \left[ |\nabla^{TZ} \eta |^2 (z)
\: - \const
\left( \parallel T_b \parallel^2
\: + \: \parallel \Pi_b \parallel^2 \right) |\eta|^2(z) \right] 
\: d\vol_{Z_b} \: d\vol_B.
\end{equation}

On a given fiber $Z_b$, for $\eta_{Z_b} \in \Omega^*(Z_b)$,
we have
\begin{equation}
\langle \eta_{Z_b}, \triangle^{Z_b} \eta_{Z_b} \rangle_{Z_b} \: = \:
\int_{Z_b} |\nabla^{TZ_b} \eta_{Z_b}|^2 (z) \: d\vol_{Z_b} \: + \: 
\sum_{ijkl} \:
\int_{Z_b} R^{Z_b}_{ijkl} \langle E^i \: I^j \: \eta, E^k \: I^l \: 
\eta \rangle
\: d\vol_{Z_b}.
\end{equation}
If $\eta_{Z_b} \in \Ker \left({\cal P} \right)$ then
\begin{equation}
\langle \eta, \triangle^{Z_b} \eta \rangle_{Z_b} \: \ge \:
A \: \diam(Z_b)^{-2} \: \langle \eta, \eta \rangle_{Z_b}.
\end{equation}
Hence
\begin{equation} \label{eq4.?}
\int_{Z_b} |\nabla^{TZ_b} \eta_{Z_b}|^2 (z) \: d\vol_{Z_b}
\: \ge \: \left(A \: \diam(Z_b)^{-2} \:  \: - \: \const \: 
\parallel R^{Z_b} \parallel_\infty \right) \langle \eta, \eta \rangle_{Z_b}
\end{equation}

From (\ref{eq4.18}), (\ref{eq4.23}) and (\ref{eq4.?}), 
if $\eta \in \Ker \left( {\cal P}^{fib} \right)$ then
\begin{equation}
\langle \eta, \triangle^M \eta \rangle_M \: \ge \:
\left( A \: \diam(Z)^{-2} \: - \: \const \left(
\parallel R^M \parallel_\infty \: + \: 
\parallel T \parallel_\infty^2 \: + \: \parallel \Pi \parallel_\infty^2 
\: + \: \parallel R^Z \parallel_\infty \right) \right) \:
\langle \eta, \eta \rangle_M.
\end{equation}
Using the Gauss-Codazzi equation, we can estimate 
$\parallel R^{Z} \parallel_\infty$
in terms of $\parallel R^M \parallel_\infty$ and 
$\parallel \Pi \parallel_\infty^2$.
The theorem follows.

\section{Collapsing to a Smooth Base} \label{sect5}

In this section we prove Theorem \ref{th2}, concerning the spectrum of the
Laplacian $\triangle^M$ on a manifold $M$ which is Gromov-Hausdorff close
to a smooth manifold $B$. We prove Theorem \ref{th3}, showing that the pairs 
$(A^\prime, h^E)$ which appear in the conclusion of Theorem \ref{th2} 
satisfy a compactness property.  
We then prove
Theorem \ref{th4}, relating the spectrum of $\triangle^M$ to the spectrum
of the differential form Laplacian on the base space $B$. \\ \\
{\bf Proof of Theorem \ref{th2} :}
For simplicity, we will omit reference to $p$. Let $g_0^{TM}$ denote
the Riemannian metric on $M$.
From \cite{Dodziuk (1982)} or Lemma \ref{lemma4}, 
if a Riemannian metric $g_1^{TM}$ on $M$ is $\epsilon$-close
to $g_0^{TM}$ then the spectrum of $\triangle^M$, computed with $g_1^{TM}$,
is $J\epsilon$-close to the spectrum computed with $g_0^{TM}$. 
We will use the geometric results of \cite{Cheeger-Fukaya-Gromov (1992)} to
find a metric $g_2^{TM}$ on $M$ which is close to $g_0^{TM}$
and to which we can apply Theorem \ref{th1}.

First, as in \cite[(2.4.1)]{Cheeger-Fukaya-Gromov (1992)}, by the
smoothing results of Abresch and others
\cite[Theorem 1.12]{Cheeger-Fukaya-Gromov (1992)},
we can find metrics
on $M$ and $B$ which are $\epsilon$-close to the original metrics such that
the new metrics satisfy $\parallel \nabla^i R \parallel_\infty \: \le \:
A_i(n, \epsilon)$ for some appropriate sequence 
$\{A_i(n, \epsilon)\}_{i=0}^\infty$. By rescaling, we may assume that
$\parallel R^M \parallel_\infty \: \le \: 1$,
$\parallel R^B \parallel_\infty \: \le \: 1$ and $\inj(B) \ge 1$.
Let $g_1^{TM}$ denote the new metric on $M$. We now
apply \cite[Theorem 2.6]{Cheeger-Fukaya-Gromov (1992)}, with $B$ fixed.
It implies that
there are positive constants $\lambda(n)$ and $c(n,\epsilon)$ so that if
$d_{GH}(M, B) \: \le \: \lambda(n)$ then there is a fibration
$f : M \rightarrow B$ such that\\
1. $\diam\left( f^{-1}(b) \right) \: \le \: c(n, \epsilon) \: d_{GH}(M, B)$.\\
2. $f$ is a $c(n, \epsilon)$-almost Riemannian submersion.\\
3. $\parallel \Pi_{f^{-1}(b)} \parallel_\infty \: \le \: 
c(n, \epsilon)$.\\
As in \cite{Fukaya (1989)}, the Gauss-Codazzi equation, 
the curvature bound on $M$ and the
second fundamental form bound on $f^{-1}(b)$ imply a
uniform bound on 
$\left\{\parallel R^{f^{-1}(b)} \parallel_\infty \right\}_{b \in B}$. 
Along with
the diameter bound on $f^{-1}(b)$, this implies that if $d_{GH}(M, B)$ is
sufficiently small then $f^{-1}(b)$ is almost flat. 

From
\cite[Propositions 3.6 and 4.9]{Cheeger-Fukaya-Gromov (1992)}, we can find
another metric $g_2^{TM}$ on $M$ which is $\epsilon$-close to $g_1^{TM}$ 
so that the fibration
$f : M \rightarrow B$ gives $M$ the structure of a Riemannian 
affine fiber bundle.
Furthermore, by \cite[Proposition 4.9]{Cheeger-Fukaya-Gromov (1992)},
there is a sequence $\{A^\prime_i(n, \epsilon)\}_{i=0}^\infty$ so that we may
assume that $g_1^{TM}$ and $g_2^{TM}$ are close in the sense that
\begin{equation} \label{eq5.1}
\parallel \nabla^i \left( g_1^{TM} - g_2^{TM} \right) \parallel_\infty
\: \le \: A^\prime_i(n, \epsilon) \: d_{GH}(M,B),
\end{equation}
where the covariant derivative in (\ref{eq5.1})
is that of the Levi-Civita connection of
$g_2^{TM}$. 
(See also \cite[Theorem 1.1]{Rong (1996)} for an explicit statement.) 
In particular, there is an upper bound on $\parallel R^M(g_2^{TM}) 
\parallel_\infty$ in terms of $B$, $n$, $\epsilon$ and $K$.

We now apply Theorem \ref{th1} to the Riemannian 
affine fiber bundle with metric
$g_2^{TM}$. It remains
to estimate the geometric terms appearing in (\ref{eq1.4}). 
We have an estimate
on $\parallel \Pi \parallel_\infty$ as above. Applying
O'Neill's formula \cite[(9.29c)]{Besse (1987)} 
to the Riemannian affine fiber bundle, 
we can estimate $\parallel T \parallel_\infty^2$ in terms
of $\parallel R^M \parallel_\infty$ and $\parallel R^B \parallel_\infty$.
Putting this together, the theorem follows. $\square$ \\

The vector bundles $E$ and Euclidean inner products $h^E$ which appear in 
Theorem \ref{th2} are not completely arbitrary.  For example, $E^0$ is
the trivial $\R$-bundle on $B$. More substantially, if $E$ is a 
real $\Z$-graded
topological vector bundle on $B$, let ${\cal C}_E$ be the space of
grading-preserving connections on $E$,
let ${\cal G}_E$ be the group of smooth grading-preserving
$\GL(E)$-gauge transformations on $E$ and
let ${\cal H}_E$ be the space of graded Euclidean inner products on 
$E$. We equip ${\cal C}_E$ and ${\cal H}_E$
with the $C^\infty$-topology. Give $({\cal C}_{E} \times {\cal H}_{E})/
{\cal G}_{E}$ the quotient topology. Let $\nabla^E$ denote the
connection part $A^\prime_{[1]}$ of $A^\prime$.

\begin{proposition} \label{prop3}
In Theorem \ref{th2}, we may assume that $E$ is one of a finite number
of isomorphism classes of
real $\Z$-graded topological vector bundles $\{E_i\}$ on $B$. 
Furthermore, there are compact subsets $C_{E_i} \subset
({\cal C}_{E_i} \times {\cal H}_{E_i})/{\cal G}_{E_i}$ 
depending on $n$, $\epsilon$ and $K$, such that we 
may assume that the gauge-equivalence class of the pair
 $\left(\nabla^E,h^E \right)$
lies in $C_E$.
\end{proposition}
\begin{pf}
As the infinitesimal holonomy of the connection $T^HM$ lies in 
$\aff(Z) = {\frak n}_R^F$, its action on ${\frak n}$, which is through the
adjoint representation, is nilpotent.
Hence its action on $\Lambda^* ( {\frak n}^* )^F$ is also
nilpotent. Given $b \in B$, it follows that the local holonomy group of $T^HM$
at $b$ acts unipotently on $E_b$.
Then the structure group of $E$ can be topologically
reduced to a discrete group and so $E$ admits a flat connection.
The dimension of $E$ is at most $2^{\dim(M) - \dim(B)}$.
By an argument of Lusztig, only a finite number of isomorphism 
classes of real topological 
vector bundles over $B$ of a given dimension
admit a flat connection \cite[p. 22]{Gromov (1982)}. This proves the
first part of the proposition.

To prove the second part of the proposition, we will first reduce to the case
$F = \{e\}$. Recall that $\widehat{Z}$ is a nilmanifold which 
covers $Z$, with covering group $F$.
Given $g \in \Aff(Z)$, we can lift it to some $\widehat{g} \in
\Aff(\widehat{Z})$. There is a automorphism $\alpha_{\widehat{g}} 
\in \Aut(F)$ such that
for all $f \in F$ and $\widehat{z} \in \widehat{Z}$, 
\begin{equation} \label{eq5.2}
\alpha_{\widehat{g}} (f) \cdot \widehat{z} 
\: = \: \left( \widehat{g} f \widehat{g}^{-1} \right)
(\widehat{z}).
\end{equation}
Considering the different possible liftings of $g$, we obtain
a well-defined homomorphism $\Aff(Z) \rightarrow \Out(F)$. Then there is
an exact sequence
\begin{equation} \label{eq5.3}
1 \longrightarrow \Aff(\widehat{Z})^F \longrightarrow \Aff(Z) \longrightarrow
\Out(F).
\end{equation}

Let $P$ be the principal $\Aff(Z)$-bundle such that $M = P \times_{Aff(Z)} Z$.
Put $\widehat{M} = P \times_{Aff(\widehat{Z})^F} \widehat{Z}$. Then 
$\widehat{M}$ is an affine fiber bundle which regularly
covers $M$, with the order of
the covering group bounded in terms of $|F|$. 
Again, there is a uniform upper bound on $|F|$ in
terms of $\dim(Z)$ \cite{Buser-Karcher (1981)}. Instead of considering $M$,
it suffices to consider
$\widehat{M}$ and work equivariantly with respect to the covering group.
Thus we assume that $Z$ is a nilmanifold, with $\Gamma \subset N_L$ and
$F = \{ e \}$.

As the fiber of $E^j$ is $\Lambda^j({\frak n}^*)$, it suffices to
prove the second part of the proposition for 
$E^1$, with fiber ${\frak n}^*$. Let us consider instead for a moment
$(E^1)^*$, with fiber ${\frak n}$. 
With respect to
the lower central series (\ref{eq1.8}) of ${\frak n}$, 
let $(E^1)^*_{[k]}$ be the vector bundle associated to $P$ with fiber
${\frak n}^\prime_{[k]}$. Then there is a filtration
\begin{equation} \label{filtration}
(E^1)^* \: = \: (E^1)^*_{[0]} \: \supset \: (E^1)^*_{[1]}
\: \supset \: \ldots \: \supset \: (E^1)^*_{[S]} \:  \supset \: 0.
\end{equation}
Let $Spl$ be the set of splittings of the short exact sequences
\begin{equation} \label{splitting}
0 \longrightarrow (E^1)^*_{[k+1]} \longrightarrow (E^1)^*_{[k]}
\longrightarrow  (E^1)^*_{[k]}/ (E^1)^*_{[k+1]} \longrightarrow 0.
\end{equation}
Put ${\cal V}_{[k]} \: = \: (E^1)^*_{[k]}/ (E^1)^*_{[k+1]}$ and 
\begin{equation} \label{V}
{\cal V} \: = \: \bigoplus_{k=0}^S {\cal V}_{[k]}.
\end{equation}
Let ${\cal H}_{\cal V}$ be the set of graded Euclidean inner products
on the $\Z$-graded real vector bundle
${\cal V}$.
A Euclidean inner product $h^{(E^1)^*}$ determines 
splittings $\{s_k\}_{k=0}^{S-1}$ of (\ref{splitting}) and a Euclidean
inner product
$h^{\cal V} \in {\cal H}_{\cal V}$. 
Conversely, one recovers $h^{(E^1)^*}$ from the
splittings $\{s_k\}_{k=0}^{S-1}$ and $h^{\cal V}$.
Thus there is an isomorphism
${\cal H}_{(E^1)^*} \cong Spl \times {\cal H}_{\cal V}$

Let ${\cal C}_{fil}$ denote the set of connections on 
$(E^1)^*$ which preserve the filtration (\ref{filtration}).
Let ${\cal C}_{{\cal V}}$ be the set of
connections on ${\cal V}$ which are
grading-preserving with respect to (\ref{V}).
Let
$\End^<({\cal V})$ be the set of endomorphisms of ${\cal V}$ which are
strictly lower-triangular with respect to $(\ref{V})$. 
Given $\left( \nabla^{(E^1)^*}, h^{(E^1)^*} \right) \in
{\cal C}_{fil} \times {\cal H}_{(E^1)^*}$, let
$i \: : \: (E^1)^* \rightarrow {\cal V}$ be the isomorphism induced
by $h^{(E^1)^*}$. Then $i \circ \nabla^{(E^1)^*} \circ i^{-1} \in 
{\cal C}_{{\cal V}} \times \Omega^1 \left(B; \End^<({\cal V}) \right)$.
In this way there is an isomorphism
\begin{equation}
{\cal C}_{fil} \times {\cal H}_{(E^1)^*} \cong
{\cal C}_{{\cal V}} \times \Omega^1 \left(B; \End^<({\cal V}) \right)
\times Spl \times {\cal H}_{\cal V}.
\end{equation}
 
Let ${\cal G}_{fil}$ be the set of filtration-preserving
gauge transformations of $(E^1)^*$ and
let ${\cal G}_{{\cal V}}$ be the set of grading-preserving
gauge transformations of ${\cal V}$.
Note that the set of splittings of (\ref{splitting}) is acted upon
freely and transitively by the gauge transformations of
$(E^1)^*_{[k]}$ which preserve
$(E^1)^*_{[k+1]}$ and act as the identity on
$(E^1)^*_{[k]}/ (E^1)^*_{[k+1]}$. It follows that
$\Ker \left( {\cal G}_{fil} \rightarrow {\cal G}_{{\cal V}} \right)$ acts
freely and transitively on $Spl$. Then
\begin{equation}
\left({\cal C}_{fil} \times {\cal H}_{(E^1)^*} \right)/
\Ker \left( {\cal G}_{fil} \rightarrow {\cal G}_{{\cal V}} \right) \cong
{\cal C}_{{\cal V}} \times \Omega^1 \left(B; \End^<({\cal V}) \right)
 \times {\cal H}_{\cal V}
\end{equation}
and so
\begin{equation} \label{todualize}
\left({\cal C}_{fil} \times {\cal H}_{(E^1)^*} \right)/{\cal G}_{fil} \cong
\left({\cal C}_{{\cal V}} \times \Omega^1 \left(B; \End^<({\cal V}) \right)
 \times {\cal H}_{\cal V} \right)/{\cal G}_{{\cal V}}.
\end{equation}
There is an obvious continuous map
$\left({\cal C}_{fil} \times {\cal H}_{(E^1)^*} \right)/{\cal G}_{fil} 
\rightarrow
\left({\cal C}_{(E^1)^*} \times {\cal H}_{(E^1)^*} \right)/{\cal G}_{(E^1)^*}$.

As $\Aff(Z)$ preserves
the lower-central-series filtration of ${\frak n}$, in our case the
dual connection to $\nabla^{E^1}$ lies in ${\cal C}_{fil}$. Then 
considering dual spaces in 
(\ref{todualize}), it is enough for us 
to show that there is a compact subset of
\begin{equation} \label{eq5.7}
\left( {\cal C}_{{\cal V}^*} \times
\Omega^1 \left( B; \End^>({\cal V}^*) \right) \times 
{\cal H}_{{\cal V}^*} \right)/{\cal G}_{{\cal V}^*}
\end{equation}
in which we may assume
that the gauge-equivalence of the pair $\left( \nabla^{E^1}, h^{E^1} \right)$
lies.  We can then map the compact subset into
$\left( {\cal C}_{E^1} \times {\cal H}_{E^1} \right)/{\cal G}_{E^1}$.

As the local holonomy of 
$E^1$ comes from an $N_R$-action, it factors through the coadjoint action 
of $N$ on ${\frak n}^*$.
Letting
$\nabla^{{\cal V}^*} = \bigoplus_{k=1}^S \nabla^{{\cal V}^*_{[k]}}$ be the
component of $\nabla^{E^1}$ in ${\cal C}_{{\cal V}^*}$, 
it follows that 
the local holonomy of $\nabla^{{\cal V}^*_{[k]}}$
is trivial and so $\nabla^{{\cal V}^*_{[k]}}$ is flat.
We first claim that there is a compact subset $C_{{\cal V}^*_{[k]}} 
\subset \left( {\cal C}_{{\cal V}^*_{[k]}} \times 
{\cal H}_{{\cal V}^*_{[k]}} \right)/
{\cal G}_{{\cal V}^*_{[k]}}$ such
that we may assume that the gauge-equivalence class of the pair
$\left( \nabla^{{\cal V}^*_{[k]}}, h^{{\cal V}^*_{[k]}} \right)$ lies in 
$C_{{\cal V}^*_{[k]}}$.

For simplicity of notation, fix $k \in [0, S]$ and let ${\cal E}$ denote
${\cal V}^*_{[k]}$. Let ${\cal F}_{\cal E}$ be the space of flat connections
on ${\cal E}$, with the subspace topology from ${\cal C}_{\cal E}$. We will
show that there is a compact subset of
$({\cal F}_{\cal E} \times 
{\cal H}_{\cal E})/{\cal G}_{\cal E}$ in which
we may assume that the gauge-equivalence class of the pair
$\left( \nabla^{\cal E}, h^{\cal E} \right)$ lies.
Then the claim will follow from mapping the compact subset into
$({\cal C}_{\cal E} \times 
{\cal H}_{\cal E})/{\cal G}_{\cal E}$.

Let $\widetilde{\cal E}$ be the lift of
${\cal E}$ to the universal cover $\widetilde{B}$ of $B$. Fix a basepoint
$\widetilde{b}_0 \in \widetilde{B}$ with projection $b_0 \in B$, 
and let $\widetilde{\cal E}_{\widetilde{b}_0}$
be the fiber of $\widetilde{\cal E}$ over $\widetilde{b}_0$. Then a flat
connection $\nabla^{\cal E}$ gives a trivialization $\widetilde{\cal E} = 
\widetilde{B} \times \widetilde{\cal E}_{\widetilde{b}_0}$. Let 
$\rho : \pi_1(B,b_0) \rightarrow \Aut(\widetilde{\cal E}_{\widetilde{b}_0})$ 
be the
holonomy of $\nabla^{\cal E}$. Then a Euclidean inner product $h^{\cal E}$ on 
${\cal E}$
can be identified with a Euclidean inner product $h^{\widetilde{\cal E}}$ on
$\widetilde{\cal E}$ which satisfies 
\begin{equation} \label{eq5.8}
h^{\widetilde{\cal E}}(\gamma^{-1} \: \widetilde{b}) \: = \: \rho(\gamma)^T \: 
h^{\widetilde{\cal E}}(\widetilde{b}) \: \rho(\gamma)
\end{equation}
for all $\gamma \in \pi_1(B, b_0)$ and $\widetilde{b} \in \widetilde{B}$.
In short, we can identify 
$({\cal F}_{\cal E} \times {\cal H}_{\cal E})/{\cal G}_{\cal E}$ with
the pairs $\left( \rho, h^{\widetilde{\cal E}} \right)$ satisfying 
(\ref{eq5.8}), 
modulo
$\Aut( \widetilde{\cal E}_{\widetilde{b}_0})$. We can use the
$\Aut( \widetilde{\cal E}_{\widetilde{b}_0})$-action to identify 
$\widetilde{\cal E}_{\widetilde{b}_0}$ with $\R^N$, with the standard 
inner product $h^{\R^N}$. If we put
\begin{align} \label{eq5.9}
X_{\cal E} \: = & \left\{ \left( \rho, h^{\widetilde{\cal E}} \right) \in 
\Hom(\pi_1(B,b_0), \GL(N, \R)) \times
{\cal H}_{\widetilde{B} \times \R^N} \: : \: 
h^{\widetilde{\cal E}}(\widetilde{b}_0) =
h^{\R^N} \text{ and for all } \right. \\
& \left.  \: \: \: \: \: \: \: \: \: \:
\gamma \in \pi_1(B, b_0) \text{ and } 
\widetilde{b} \in \widetilde{B}, \: 
h^{\widetilde{\cal E}}(\gamma^{-1} \:  \widetilde{b}) \: = \: 
\rho(\gamma)^T \: 
h^{\widetilde{\cal E}}(\widetilde{b}) \: \rho(\gamma)
\right\} \notag
\end{align}
then we have identified $({\cal F}_{\cal E} \times 
{\cal H}_{\cal E})/{\cal G}_{\cal E}$ with
$X_{\cal E}/O(N)$.
Let $\{\gamma_j\}$ be a finite generating set of $\pi_1(B,b_0)$.
The topology on $X_{\cal E}$ comes from the fiber bundle structure 
\begin{equation} \label{eq5.10}
X_{\cal E} \rightarrow \Hom(\pi_1(B,b_0), \GL(N, \R)),
\end{equation}
whose fiber over $\rho \in
\Hom(\pi_1(B,b_0), \GL(N, \R))$ is 
\begin{align} \label{eq5.11}
& \left\{ h^{\widetilde{\cal E}} \in 
{\cal H}_{\widetilde{B} \times \R^N} \: : 
h^{\widetilde{\cal E}}(\widetilde{b}_0) =
h^{\R^N} \text{ and for all $\gamma \in \pi_1(B, b_0)$ and 
$\widetilde{b} \in \widetilde{B}$,} \right. \\
& \left. \: \: \: \: \: \: \: \: \: \:
h^{\widetilde{\cal E}}(\gamma^{-1} \:  \widetilde{b}) \: = \: 
\rho(\gamma)^T \: 
h^{\widetilde{\cal E}}(\widetilde{b}) \: \rho(\gamma)
\right\} \notag
\end{align}
Here
$\Hom(\pi_1(B,b_0), \GL(N, \R))$ has a topology as a subspace of
$\GL(N, \R)^{\{\gamma_j\}}$
and the fiber (\ref{eq5.11}) has the $C^\infty$-topology.
Thus it suffices to show that
$\left( \rho, h^{\widetilde{\cal E}} \right)$ lies in a predetermined compact
subset $C_{\cal E}$ of $X_{\cal E}$.

By \cite[(1-7)]{Fukaya (1989)}, we may assume that we have uniform bounds
on the second fundamental form $\Pi$ of the Riemannian
affine fiber bundle $M$, along with
its covariant derivatives. As $\Pi$ determines how the Riemannian
metrics on nearby
fibers vary (with respect to $T^HM$), 
and $h^{E^1}$ comes from the inner product on
the parallel differential forms on the fibers $\{Z_b\}_{b \in B}$, we 
obtain uniform bounds on $\left( h^{E_1} \right)^{-1} \left( 
\nabla^{E_1} h^{E_1}
\right)$ and its covariant derivatives.
In particular, we also have a uniform bound
on $\left( h^{\cal E} \right)^{-1} \left( 
\nabla^{\cal E} h^{\cal E}
\right)$ and hence on
$\left( h^{\widetilde{\cal E}} \right)^{-1} \: \left( d \: 
h^{\widetilde{\cal E}}
\right)$. 
For the 
finite generating set $\{\gamma_j\}$, 
using the fact that 
$h^{\widetilde{\cal E}}(\widetilde{b}_0) =
h^{\R^N}$,
we obtain in this way uniform bounds
on $\{ h^{\widetilde{\cal E}} (\gamma_j^{-1} \: \widetilde{b}_0) \}$.
The equivariance
(\ref{eq5.8}) then gives uniform bounds on $\{ \rho(\gamma_j)^T \: 
\rho(\gamma_j)\}$ and hence on $\{\rho(\gamma_j)\}$. Thus 
$\rho$ lies in a predetermined compact subset of the representation space
$\Hom(\pi_1(B,b_0), \GL(N, \R))$. Given $\rho$, 
the
uniform bounds on the covariant derivatives of $h^{\widetilde{\cal E}}$ over
a fundamental domain in $\widetilde{B}$ show that we have compactness in
the fiber (\ref{eq5.11}). 
As these bounds can be made continuous in $\rho$, the claim 
follows.

Fix a Euclidean inner product $h^{{\cal V}^*}_0$ on ${\cal V}^*$. Given a pair
$\left( \nabla^{{\cal V}^*}, h^{{\cal V}^*} \right) \in {\cal C}_{{\cal V}^*} 
\times {\cal H}_{{\cal V}^*}$, 
we can always perform a gauge transformation to transform
the Euclidean inner product 
to $h^{{\cal V}^*}_0$. Let ${\cal O}_{{\cal V}^*}$ be 
the orthogonal
gauge transformations with respect to $h^{{\cal V}^*}_0$. Then we can identify
$\left( {\cal C}_{{\cal V}^*} \times {\cal H}_{{\cal V}^*} \right)/
{\cal G}_{{\cal V}^*}$ with
${\cal C}_{{\cal V}^*} / {\cal O}_{{\cal V}^*}$. Similarly,
\begin{equation} \label{eq5.12}
\left( {\cal C}_{{\cal V}^*} \times
\Omega^1 \left( B; \End^>({\cal V}^*) \right) \times 
{\cal H}_{{\cal V}^*} \right)/{\cal G}_{{\cal V}^*} \cong
\left( {\cal C}_{{\cal V}^*} \times
\Omega^1 \left( B; \End^>({\cal V}^*) \right)
\right)/{\cal O}_{{\cal V}^*}.
\end{equation}
There is a singular fibration
$p : \left( {\cal C}_{{\cal V}^*} \times
\Omega^1 \left( B; \End^>({\cal V}^*) \right)
\right)/{\cal O}_{{\cal V}^*} \rightarrow
{\cal C}_{{\cal V}^*} 
/{\cal O}_{{\cal V}^*}$. The fiber over a gauge-equivalence class
$[\nabla^{{\cal V}^*}]$ is 
$\Omega^1 \left( B; \End^>({\cal V}^*) \right)/G$, where $G$ is the
group of orthogonal gauge transformations which are parallel with
respect to $\nabla^{{\cal V}^*}$. In particular, upon choosing a basepoint 
$b_0 \in B$, we can view
$G$ as contained in the finite-dimensional orthogonal group
$O({\cal V}^*_{b_0})$.

From what we have already shown, 
we know that we are restricted to a compact subset of the base 
$\left( {\cal C}_{{\cal V}^*} \times {\cal H}_{{\cal V}^*} \right)/
{\cal G}_{{\cal V}^*} \cong {\cal C}_{{\cal V}^*} 
/{\cal O}_{{\cal V}^*}$ of the singular fibration $p$.
Let $\left( \nabla^{E^1}
\right)^T$ be the adjoint connection to $\nabla^{E^1}$ with
respect to $h^{E^1}$. 
The uniform bounds on $\left( h^{E^1} \right)^{-1} \left( 
\nabla^{E^1} h^{E^1}
\right)$ and its derivatives give uniform $C^\infty$-bounds on the
part of $\nabla^{E^1}$ which does not preserve the metric
$h^{E^1}$, i.e. on  
$\nabla^{E^1} - \left( \nabla^{E^1}
\right)^T \in \Omega^1(B; \End(E^1))$. In particular,
using the upper triangularity of $\nabla^{E^1}$,
we obtain uniform $C^\infty$-bounds on the part of $\nabla^{E^1}$
in $\Omega^1 \left( B; \End^>({\cal V}^*) \right)$.  As the bounds 
can be made continous with respect to $\left[ \nabla^{{\cal V}^*} \right] \in
{\cal C}_{{\cal V}^*} /{\cal O}_{{\cal V}^*}$, we have shown that there is
a fixed compact subset of
$\left( {\cal C}_{{\cal V}^*} \times
\Omega^1 \left( B; \End^>({\cal V}^*) \right) \times 
{\cal H}_{{\cal V}^*} \right)/{\cal G}_{{\cal V}^*}$ in which we may assume
that the pair $\left( \nabla^{E^1}, h^{E^1} \right)$ lies.

To summarize, we have shown that the topological vector bundle $E^1$ has
a flat structure ${\cal V}^*$, with flat connection $\nabla^{{\cal V}^*}$.
We showed that there are bounds on the holonomy of $\nabla^{{\cal V}^*}$  
which are uniform in $n$, $\epsilon$ and $K$.  We then showed that
$h^{E^1}$ and $\nabla^{E^1} - \nabla^{{\cal V}^*}$ are $C^\infty$-bounded
in terms of $n$, $\epsilon$ and $K$. (More precisely, we showed that these
statements are true after an appropriate gauge transformation is made.)
The proposition follows.
\end{pf}

Let ${\cal S}_{E}$ be the space of degree-$1$
superconnections on $E$, with the $C^\infty$-topology.

\begin{proposition} \label{prop4}
With reference to Proposition \ref{prop3},
there are compact subsets $D_{E_i} \subset
({\cal S}_{E_i} \times {\cal H}_{E_i})/{\cal G}_{E_i}$ 
depending on $n$, $\epsilon$ and $K$ such that we 
may assume that the gauge-equivalence class of the pair
 $\left(A^\prime,h^E \right)$
lies in $D_E$.
\end{proposition}
\begin{pf}
Let $E$ be as in Proposition \ref{prop3}.
As in the proof of Proposition \ref{prop3}, 
upon choosing $h^E_0$, we have identifications 
$({\cal C}_E \times {\cal H}_E)/{\cal G}_E \cong 
{\cal C}_E/{\cal O}_E$ and
$({\cal S}_E \times {\cal H}_E)/{\cal G}_E \cong 
{\cal S}_E/{\cal O}_E$. There is a singular fibration 
$p : {\cal S}_E/{\cal O}_E \rightarrow {\cal C}_E/{\cal O}_E$ coming
from the projection $A^\prime \rightarrow A^\prime_{[1]}$. The fiber of $p$
over a gauge equivalence class $\left[ \nabla^E \right]$ is
$\bigoplus_{a+b=1} \Omega^a(B; \End(E^*, E^{*+b}))/G$, where
$G$ is the group of orthogonal gauge transformations which are parallel
with respect to $\nabla^E$.

From Proposition \ref{prop3}, we know that we are restricted to a
compact subset of the base $({\cal C}_E \times {\cal H}_E)/{\cal G}_E \cong 
{\cal C}_E/{\cal O}_E$. 
The superconnection on $E$ has the
form (\ref{eq4.10}). We measure norms on $\Omega(B; \End(E))$ using $h_0^E$.
The differential $d^{\frak n}$ comes from exterior 
differentiation on the parallel forms on the fibers of the Riemannian
affine fiber
bundle. Note that as $A^\prime$ is flat,
$d^{\frak n}$ is parallel with respect to $\nabla^E$.  
As we have a uniform ($n$, $\epsilon$ and $K$)-dependent bound on the
curvatures of the fibers $Z$,
Lemma \ref{lemma2} gives
a uniform bound on the structure constants $\{ c^i_{\: jk}\}$ and hence
a uniform bound on $\parallel d^{\frak n} \parallel_\infty$.
The operator $i_T$ is also parallel with respect to $\nabla^E$. 
From O'Neill's formula \cite[(9.29)]{Besse (1987)}, we 
obtain a uniform bound on $\parallel i_T \parallel_\infty$.
Thus we have uniform $C^\infty$-bounds on $A^\prime - \nabla^E \in 
\Omega(B; \End(E))$ and so we have compactness in the fibers of $p$.
As the bounds can be made continuous
with respect to $\left[ \nabla^E \right] \in
{\cal C}_E/{\cal O}_E$, the proposition follows.
\end{pf}

Propositions \ref{prop3} and \ref{prop4} prove Theorem \ref{th3}.

We will need certain eigenvalue statements.
Let ${\cal E}$ be a $\Z_2$-graded real topological
vector bundle on a smooth closed manifold $B$.
Let ${\cal S}_{\cal E}$ be the space of
superconnections on ${\cal E}$,
let ${\cal G}_{\cal E}$ be the $\GL({\cal E})$-gauge group of ${\cal E}$ and
${\cal H}_{\cal E}$ be the space of Euclidean metrics on ${\cal E}$. 
Fix a Euclidean metric $h^{\cal E}_0 \in {\cal H}_{\cal E}$.
Given a pair $\left( A^\prime, h^{\cal E} \right) \in 
{\cal S}_{\cal E} \times {\cal H}_{\cal E}$,
we can always perform a gauge transformation to transform the Euclidean metric
to $h^{\cal E}_0$. Let ${\cal O}_{\cal E}$ be the group of orthogonal gauge
transformations of ${\cal E}$ with respect to $h^{\cal E}_0$.
Then we can identify
$({\cal S}_{\cal E} \times {\cal H}_{\cal E})/{\cal G}_{\cal E}$ with 
${\cal S}_{\cal E}/{\cal O}_{\cal E}$. 
Given $A^\prime \in {\cal S}_{\cal E}$, let
$\left (A^\prime \right)^*$ be its adjoint with respect to
$h^{\cal E}_0$ and put 
$\triangle_{A^\prime} \: = \: A^\prime \: \left (A^\prime \right)^* \: + \: 
\left (A^\prime \right)^* \: A^\prime$, acting on
$\Omega(B; {\cal E})$. 
For $j \in \Z^+$, let $\mu_j \left(A^\prime \right)$ be the 
$j$-th eigenvalue of $\triangle_{A^\prime}$, 
counted with multiplicity. It is ${\cal O}_{\cal E}$-invariant.
Equivalently, $\mu_j$ is ${\cal G}_{\cal E}$-invariant
as a function of the pair $(A^\prime, h^E)$. \\ \\
{\bf Proof of Theorem \ref{th4} :}
As $E^p$ admits a flat connection
there is some $r \in \N$ such that for all $p$, $E^p \otimes \R^r$ is
topologically isomorphic to the trivial vector bundle
$B \times \R^{r \cdot rk(E^p)}$. Hence
$E \otimes \R^r$ is topologically isomorphic to the
$\Z$-graded trivial vector bundle 
${\cal E} = B \times \R^{r \cdot rk(E)}$.

For simplicity, we omit reference to $p$. In view of Theorem \ref{th2}, it
suffices to show that there is a positive constant $D(n,\epsilon, K)$ such
that $|\lambda_j(B; E)^{1/2} \: - \:
\lambda^\prime_j(B)^{1/2}| \: \le \: D(\epsilon, n, K)$.

The operator $\triangle^E \otimes \Id$ on 
$\Omega(B; E) \otimes \R^r$ has a spectrum which is the same as that
of $\triangle^E$, but with multiplicities multiplied by $r$.
Hence it is enough to compare the spectrum of $\triangle^E \otimes \Id$,
acting on $\Omega(B; E \otimes \R^r)$, with that of the standard Laplacian
on $\Omega(B; {\cal E})$.

From Theorem \ref{th3}, we may assume that the 
gauge-equivalence class of the pair $(A^\prime, h^E)$ lies in
a predetermined compact subset $D \subset 
({\cal S}_{E} \times {\cal H}_{E})/
{\cal G}_{E}$. Put $A_1^\prime = A^\prime \otimes \Id$, 
acting on $\Omega(B; E \otimes \R^r)$, and put $h_1 = h^E \otimes
h^{\R^r}$ on $E \otimes \R^r$. Using the isomorphism
${\cal E} \cong E \otimes \R^r$, we may assume that
the gauge-equivalence class of the pair $(A_1^\prime, h_1)$ lies in a
predetermined compact subset $D_1 \subset 
({\cal S}_{\cal E} \times {\cal H}_{\cal E})/
{\cal G}_{\cal E}$.

Let $A^\prime_2$ be the trivial flat connection
on ${\cal E}$ and let $h_2$ be the product Euclidean inner product
on ${\cal E}$.
With an appropriate gauge transformation $g \in {\cal G}_{\cal E}$, we can
transform $\left( A_1^\prime,  h_1 \right)$ to
$\left( g \cdot A^\prime_1, 
h_2 \right)$ without changing the eigenvalues.
Under the identification $({\cal S}_{\cal E} \times {\cal H}_{\cal E})/
{\cal G}_{\cal E} = {\cal S}_{\cal E}/{\cal O}_{\cal E}$,
we can assume that the equivalence class of $g \cdot A^\prime_1$ 
in ${\cal S}_{\cal E}/{\cal O}_{\cal E}$ lies in
a predetermined compact subset $D_2 \subset 
{\cal S}_{\cal E}/{\cal O}_{\cal E}$.

The eigenvalues of the Laplacian associated to the superconnection
$A^\prime_1$ and the Euclidean inner product $h_2$ are unchanged when
the group of orthogonal gauge transformations  ${\cal O}_{\cal E}$ acts
on $A^\prime_1$.
Consider the function $l : {\cal S}_{\cal E} \times {\cal S}_{\cal E}
\rightarrow \R$ given by
\begin{equation} \label{eq5.25}
l(A^\prime_1, A^\prime_2) \: = \: \inf_{g^\prime \in 
{\cal O}_{\cal E}} \parallel g^\prime \cdot A^\prime_1 -
A^\prime_2 \parallel.
\end{equation}
An elementary argument shows that $l$ is continuous. Hence it
descends to a continuous function on 
$\left( {\cal S}_{\cal E}/{\cal O}_{\cal E} \right) \times
\left( {\cal S}_{\cal E}/{\cal O}_{\cal E} \right)$.
Applying Lemma \ref{lemma?} to $g^\prime \cdot 
(g \cdot A^\prime_1)$ and $A^\prime_2$, 
the compactness of $D_2$ and the finiteness statement in 
Theorem 3
give the desired eigenvalue estimate.
The theorem follows. $\square$

\section{Small Positive Eigenvalues} \label{sect8}

In this section we characterize the manifolds $M$ for which 
the $p$-form Laplacian has small positive eigenvalues. We first
describe a spectral sequence which computes the cohomology of a flat
degree-$1$ superconnection $A^\prime$. We use the compactness result of
Theorem \ref{th3} to show that if $M$ has $j$ small eigenvalues of
the $p$-form Laplacian, with $j > \bb_p(M)$, and $M$ collapses to a smooth
manifold $B$ then
there is an associated flat degree-$1$ superconnection $A^\prime_\infty$ on $B$
with $\dim(\HH^p(A^\prime_\infty)) \ge j$. We then use the spectral sequence of
$A^\prime_\infty$ to characterize when this can happen.
In Corollary \ref{cor2}
we give a bound on the number of small eigenvalues
of the $1$-form Laplacian. 
In Corollary \ref{cor4} we give a bound on the number of
small eigenvalues of the $p$-form Laplacian when one is sufficiently close
to a smooth limit space of dimension $\dim(M)-1$.
Theorem \ref{th6} describes when a collapsing sequence can have small
positive
eigenvalues of the $p$-form Laplacian, in terms of the topology of the
affine fiber bundle $M \rightarrow B$.
Corollary \ref{cor5} gives a precise description of when there are small
positive eigenvalues of the $1$-form Laplacian in a collapsing sequence. 
In Corollaries \ref{cor6} and \ref{cor7}
we look at collapsing sequences with limit spaces of
dimension $1$ or $\dim(M) - 1$, respectively.
Finally, given an affine 
fiber bundle, in
Theorem \ref{th7} we give a collapsing construction which produces
small eigenvalues of the $p$-form Laplacian.

In the collapsing arguments in this section, when the limit space is a
smooth manifold, we can always assume that its Riemannian metric is smooth.
At first sight the
smoothness assumption on the metric may seem strange, as the limit space of a
bounded-sectional-curvature collapse, when a smooth manifold, generally
only has a  $C^{1,\alpha}$-metric.
The point is that we are interested in the case when an eigenvalue goes
to zero, which gives a zero-eigenvalue of $\triangle^E$ in the limit.
The property of having a zero-eigenvalue is essentially topological in
nature and so will also be true for a smoothed metric. For this reason,
we can apply smoothing results to the metrics and so ensure that the
limit metric is smooth.

Let $B$ be a smooth connected closed manifold.
Let $E = \oplus_{j=0}^m E^j$ be a 
$\Z$-graded real vector bundle on $B$ and
let $A^\prime = \sum_{i \ge 0} A^\prime_{[i]}$ be a flat
degree-$1$ superconnection on $E$.
Let $\HH^p(A^\prime)$ denote the degree-$p$
cohomology of the differential $A^\prime$
on $\Omega(B; E)$, where the latter has the total grading.
Given $a, b \in \N$, we will write
$\omega^{a,b}$ for an element of $\Omega^a(B; E^b)$.

In order to compute $\HH^p(A^\prime)$, let us first consider the equation 
$A^\prime \omega = 0$. Putting 
\begin{equation} \label{eq8.1}
\omega \: = \: \omega^{p,0} \: + \: 
\omega^{p-1,1} \: + \: \omega^{p-2,2} \: + \: \ldots,
\end{equation} 
we obtain
\begin{equation} \label{eq8.2}
\left( A^\prime_{[0]} \: + \: A^\prime_{[1]} \: + \: A^\prime_{[2]} \: + \:
\ldots \right) \left( \omega^{p,0} \: + \: 
\omega^{p-1,1} \: + \: \omega^{p-2,2} \: + \: \ldots \right) \: = \: 0,
\end{equation}
or
\begin{align} \label{eq8.3}
A^\prime_{[0]} \: \omega^{p,0} \: & = 0, \\
A^\prime_{[0]} \: \omega^{p-1,1} \: + \: A^\prime_{[1]} \: \omega^{p,0}
 & = 0, \notag \\
A^\prime_{[0]} \: \omega^{p-2,2} \: + \: A^\prime_{[1]} \: \omega^{p-1,1}
 \: + \: A^\prime_{[2]} \: \omega^{p,0}
 & = 0, \notag \\
\vdots \notag
\end{align}
We can try to solve these equations iteratively.

Formalizing
this procedure, we obtain a spectral sequence to compute $\HH^p(A^\prime)$. 
Put $E_0^{a,b} = \Omega^a(B; E^b)$ and
define $d_0 : E_0^{a,b} \rightarrow E_0^{a,b+1}$ by
$d_0 \: \omega^{a,b} = A_{[0]} \: \omega^{a,b}$.
For $r \ge 1$, put
\begin{equation} \label{eq8.4}
E_r^{a,b} = \frac{
\{ \{\omega^{a+s, b-s}\}_{s=0}^{r-1} : 
\text{ for $0 \le s \le r-1$, }
\sum_{t = 0}^s A^\prime_{[s-t]} \: \omega^{a+t, b-t} 
= 0 \}
}{
\{ \{\omega^{a+s, b-s}\}_{s=0}^{r-1} : 
\omega^{a+s,b-s} = \sum_{t = 0}^s A^\prime_{[s-t]} \: 
\widehat{\omega}^{a+t, b-t-1} \text{ for some }
\{\widehat{\omega}^{a+s, b-s-1}\}_{s=0}^{r-1}
\}
}. 
\end{equation}
Define a differential $d_r : E_r^{a,b} \rightarrow E_r^{a+r,b-r+1}$ by
\begin{equation} \label{eq8.5}
d_r \: \left\{ \omega^{a+s, b-s} \right\}_{s=0}^{r-1} \: = \: 
\left\{ \sum_{t=0}^{r-1} A^\prime_{[r+s-t]} \:
\omega^{a+t,b-t} \right\}_{s=0}^{r-1}.
\end{equation}
Then $E_{r+1} \cong \Ker(d_r)/\Image(d_r)$. The spectral sequence
$\{E_r^{*,*}\}_{r=0}^\infty$ has a limit $E_\infty^{*,*}$ with
\begin{equation} \label{eq8.6}
\HH^p(A^\prime) \cong \bigoplus_{a+b=p} E_\infty^{a,b}.
\end{equation}
From \cite[Proposition 2.5]{Bismut-Lott (1995)}, for each $b \in \N$,
$\HH^b(A_{[0]}^\prime)$ is a flat vector bundle on $B$. 
Then
\begin{align} \label{eq8.7}
E_0^{a,b} & = \Omega^a(B; E^b), \\
E_1^{a,b} & = \Omega^a(B; \HH^b(A_{[0]}^\prime)), \notag \\
E_2^{a,b} & = \HH^a(B; \HH^b(A_{[0]}^\prime))). \notag
\end{align}
\noindent
{\bf Example 4 : } If $M \rightarrow B$ is a fiber bundle, $E$ is the
infinite-dimensional vector bundle $W$ of vertical differential forms
\cite[Section III(a)]{Bismut-Lott (1995)} and $A^\prime$ is the
superconnection arising from
exterior differentiation on $M$ then we recover the Leray
spectral sequence to compute $\HH^*(M; \R)$.\\ \\
{\bf Example 5 : } If $M \rightarrow B$ is an affine fiber bundle,
$E$ is the vector bundle of parallel differential forms on the fibers 
and $A^\prime$ is as in (\ref{eq4.10}) then it follows from 
\cite[Corollary 7.28]{Raghunathan (1972)} that 
$E_r^{*,*}$ is the same as the corresponding term in the Leray
spectral sequence for $\HH^*(M; \R)$ if $r \ge 1$.\\

Suppose that $M$ is a connected closed manifold with at least $j$ small
eigenvalues of $\triangle_p$ for some
$j > \bb_p(M)$. Consider a
sequence of Riemannian metrics $\{g_i\}_{i=1}^\infty$ in ${\cal M}(M, K)$ 
with $\lim_{i \rightarrow\infty} 
\lambda_{p,j}(M, g_i) = 0$.
As in the proof of Theorem \ref{th2}, for any $\epsilon > 0$
there is a sequence $\{A_k(n, \epsilon)\}_{k=0}^\infty$ so that for
all $i$, we can
find a new metric $g^\prime_i$ on $M$ which is $\epsilon$-close to $g_i$,
with 
$\parallel \nabla^k R^{M}(g^\prime_i) \parallel_\infty \: \le \:
A_k(n, \epsilon)$. 
Fix $\epsilon$ to be, say, $\frac{1}{2}$.
From \cite{Dodziuk (1982)} or Lemma \ref{lemma4}, 
we have that $\lambda_{p,j}(M, g^\prime_i)$ is $J \epsilon$-close to
$\lambda_{p,j}(M, g_i)$ for some fixed integer $J$.
Thus without loss of generality, we may replace $g_i$ by $g^\prime_i$. We
relabel $g^\prime_i$ as $g_i$.

As $j \: > \: \bb_p(M)$, there must be a subsequence of
$\{(M, g_i)\}_{i=1}^\infty$
which Gromov-Hausdorff
converges to a lower-dimensional limit space $X$. That is, we are in 
the collapsing situation. 
Suppose that the limit space is
a smooth manifold $B$. From 
\cite[Section 5]{Cheeger-Fukaya-Gromov (1992)}, 
the regularity of the metrics on $M$
implies that $B$ has a smooth Riemannian metric
$g^{TB}$. (We are in the situation in which the limit space $\check{X}$ of
the frame bundles, a smooth Riemannian
manifold, has an $O(n)$-action with a single
orbit type.)
From Theorem \ref{th2},
for large $i$ there are vector bundles $E_i$ on $B$, flat degree-$1$
superconnections $A^\prime_i$ on $E_i$, and Euclidean inner products
$h^{E_i}$ on $E_i$ such that $\lambda_{p,j}(M, g_i)$ is $\epsilon$-close
to $\lambda_{p,j}(B; E_i)$. From Theorem \ref{th3}, after taking a
subsequence we may assume that all of the $E_i$'s are topologically
equivalent to a single vector bundle $E$ on $B$, and that
the pairs $\left( A^\prime_i, h^{E_i} \right)$ 
converge after gauge transformation to
a pair $\left( A^\prime_\infty, h^{E_\infty} \right)$. 
Then from Lemmas \ref{lemma4} and \ref{lemma?}, the Laplacian associated
to $\left( A^\prime_\infty, h^{E_\infty} \right)$ satisfies
$\dim \Ker \left( \triangle^E_p \right) \: \ge \: j$. Applying standard
Hodge theory to the superconnection Laplacian $\triangle^E$,
we obtain $\dim \left( \HH^p(A^\prime_\infty) \right) \: \ge \: j$.
On the other hand, looking at the $E_2$-term of the spectral sequence gives
$\dim \left( \HH^p(A^\prime_\infty) \right) \: \le \: 
\sum_{a+b=p} \dim \left( \HH^a(B; \HH^b(A_{\infty,[0]}^\prime) \right)$.
Thus 
\begin{equation} \label{whoknows}
j \: \le \: \sum_{a+b=p} \dim \left( \HH^a(B; \HH^b(A_{\infty,[0]}^\prime)) 
\right).
\end{equation}
{\bf Proof of Corollary \ref{cor2} : }
In the case $p \: = \: 1$, we obtain
\begin{equation}
j \: \le \: \dim \left( \HH^1(B; \HH^0(A^\prime_{\infty,[0]}) \right) +
\dim \left( \HH^0(B; \HH^1(A^\prime_{\infty,[0]}) \right).
\end{equation}
As $\HH^0(A^\prime_{\infty,[0]})$ is the trivial $\R$-bundle on $B$,
$\dim \left( \HH^1(B; \HH^0(A^\prime_{\infty,[0]}) \right) = \bb_1(B)$.
As $A^\prime_{\infty,[0]}$ acts by zero on $E^0$, 
there is an injection $\HH^1(A^\prime_{\infty,[0]}) \rightarrow E^1$. 
Then 
\begin{equation} \label{eq8.10}
\dim \left( \HH^0(B; \HH^1(A^\prime_{\infty,[0]})) \right) \: \le \:
\dim \left( \HH^1(A^\prime_{\infty,[0]}) \right) 
\le \: \dim \left( E^1 \right)
\: \le \: \dim(M) \: - \: \dim(B).
\end{equation}
Thus $j \: \le \: \bb_1(B) \: + \dim(M) \: - \: \dim(B)$.
On the other hand, the spectral sequence for $\HH^*(M; \R)$ gives
\begin{equation} \label{eq8.11}
\HH^1(M; \R) \: = \: \HH^1(B; \R) \: \oplus \: \Ker \left( 
\HH^0(B; \HH^1(Z; \R)) \rightarrow \HH^2(B; \R) \right).
\end{equation}
In particular, $\bb_1(B) \: \le \: \bb_1(M)$.
The corollary follows. $\square$ \\ \\
{\bf Remark : } Using heat equation methods \cite{Berard (1988)} one can
show that there is an increasing function $f$ such that if $\Ric(M) \: \ge \:
- (n-1) \: \lambda^2$ and $\diam(M) \: \le \: D$ then the number of small 
eigenvalues of the $1$-form Laplacian is bounded above by 
$f(\lambda \: D)$. This result is weaker than Corollary \ref{cor2} when
applied to manifolds with sectional curvature bounds, but is more general
in that it applies to manifolds with just a lower Ricci curvature bound.\\ \\
{\bf Proof of Corollary \ref{cor4} :}
From Fukaya's fibration theorem, if  
a manifold $M^n$ with $\parallel R^M \parallel_\infty \: \le \: K$ 
is sufficiently 
Gromov-Hausdorff close to $B$ then $M$ is the total space of a circle bundle
over $B$. 
Suppose that the claim of the corollary is not true.
Then there is a sequence of connected closed $n$-dimensional Riemannian 
manifolds $\{(M_i, g_i)\}_{i=1}^\infty$ with $\parallel R^{M_i}(g_i) 
\parallel_\infty \: \le \: K$ and $\lim_{i \rightarrow \infty} M_i = B$
which provides a counterexample.
As there is a finite number of isomorphism
classes of flat real line
bundles on $B$, after passing to a subsequence we may assume that
each $M_i$ is a circle bundle over $B$ with a fixed 
orientation bundle
${\cal O}$ and that $\lim_{i \rightarrow \infty} \lambda_{p,j}(M_i, g_i) = 0$
for $j \: = \: \bb_p(B) \: + \: \bb_{p-1}(B; {\cal O}) \: + \: 1$. 
Following the argument before the proof of Corollary \ref{cor2},
we obtain $E \: = \: E^0 \: \oplus \: E^1$ on $B$,
with $E^0$ a trivial $\R$-bundle and $E^1 \: = \: {\cal O}$, and a limit
superconnection $A^\prime_\infty$ on $E$ with $A^\prime_{\infty,[0]} = 0$ and
$A^\prime_{\infty,[1]} = \nabla^E$, the canonical flat connection. Then as in
(\ref{whoknows}), we obtain
\begin{equation}
j \: \le \: 
\bb_p(B) \: + \: \bb_{p-1}(B; {\cal O}),
\end{equation}
which is a contradiction. \\ \\
{\bf Proof of Theorem \ref{th6} :}
As in the proof of Theorem \ref{th2},
without loss of generality we may assume that each $(M, g_i)$ is a
Riemannian affine fiber bundle structure on the affine fiber bundle
$M \rightarrow B$.
Suppose that for each $q \in [0, p]$,
$b_q(Z) \: = \: \dim \left( \Lambda^q({\frak n}^*)^F \right)$ and the holonomy
representation of the flat vector bundle $\HH^q(Z; \R)$ on $B$ is
semisimple. 
Let $E \rightarrow B$ be the real vector bundle associated to the
affine fiber bundle $M \rightarrow B$ as in Section
\ref{sect4}. Then
$E \cong \HH^*(Z; \R)$. The superconnection $A^\prime_E$ on $E$,
from Section \ref{sect4}, has
$A^\prime_{E,[0]} = 0$ and $A^\prime_{E,[1]} = \nabla^E$, the
canonical flat connection on $E \cong \HH^*(Z; \R)$.
As the affine fiber bundle is fixed, each $E_i$ equals $E$ and each
$A^\prime_i$ equals $A^\prime_E$. However,
the Euclidean metrics $\{ h^E_i\}_{i=1}^\infty$ on $E$ vary. 
There is a sequence of gauge transformations
$\{ g_i\}_{i=1}^\infty$ so that after passing to a subsequence,
$\lim_{i \rightarrow \infty} g_i \cdot (A^\prime_i, h^E_i) =
(A^\prime_\infty, h^E_\infty)$ for some pair $(A^\prime_\infty, h^E_\infty)$. 
Clearly $A^\prime_{\infty, [0]} = 0$ and
$A^\prime_{\infty, [1]} = \lim_{i \rightarrow \infty} g_i \cdot \nabla^E$.
As the holonomy representation of $\HH^q(Z; \R)$ is semisimple
for $q \in [0, p]$, the connection
$A^\prime_{\infty, [1]} \big|_{E^q}$ is gauge-equivalent to 
$\nabla^{E^q}$. That is, the connection does not degenerate. 
(In the complex case this follows from
\cite[Theorem 1.27]{Lubotzky-Magid (1985)} and the real case follows from
\cite[Theorem 11.4]{Richardson (1988)}.) 

Equation (\ref{whoknows}) now implies that
\begin{equation}
j \: \le \: \sum_{a+b=p} \dim \left( \HH^a(B; \HH^b(Z; \R)) \right).
\end{equation}
If the Leray spectral sequence to compute $\HH^*(M; \R)$ 
degenerates at the $E_2$ term then
\begin{equation}
\bb_p(M) \: = \: \sum_{a+b=p} \dim \left( \HH^a(B; \HH^b(Z; \R)) \right),
\end{equation}
which contradicts the assumption that $j \: > \: \bb_p(M)$. $\square$ \\ \\
{\bf Example 6 : } Let $Z$ be an almost flat manifold as in Example 1.  Put
$M = Z \times B$. If there is a sequence of affine-parallel metrics on $Z$
which give it $r_p$ small eigenvalues of the $p$-form Laplacian
then $M$ has $\sum_{a + b =p} r_a \cdot \bb_b(B)$ small eigenvalues of the
$p$-form Laplacian.  This gives an example of Theorem \ref{th6}.1.\\ \\
{\bf Example 7 \cite{Jammes (2000)} : } Let $N$ be the Heisenberg group
of upper-diagonal unipotent $3 \times 3$ matrices and let $\Gamma$ be the
integer lattice in $G$. Put $M = \Gamma \backslash N$. Then $M$ fibers
over $S^1$, the fiber being $T^2$ and the monodromy being given by the matrix
$\begin{pmatrix}
1 & 1 \\
0 & 1
\end{pmatrix}$. One has $\bb_1(M) = 2$, but for any $K > 0$, $a_{1,3,K} = 0$.
That is, one can collapse $M$ to a circle by a sequence of affine-parallel
metrics, while producing $3$ small eigenvalues of the $1$-form Laplacian.
This gives an example of Theorem \ref{th6}.2.\\ \\
{\bf Example 8 : } Consider $M$ as in Example 2.  If the Leray spectral
sequence to compute $\HH^p(M; \R)$ does not degenerate at the $E_2$ term
then there are small positive eigenvalues of the $p$-form Laplacian on $M$.
This gives an example of Theorem \ref{th6}.3. \\ \\
{\bf Proof of Corollary \ref{cor5} :}
The affine fiber bundle $M \rightarrow B$ induces a vector bundle
$E \rightarrow B$ and a flat degree-$1$ superconnection $A^\prime_E$,
as in Section \ref{sect4}. As in Example 5, 
the spectral sequence associated to
$A^\prime_E$ is the same as the Leray spectral sequence to compute
$\HH^*(M; \R)$. Let $A^\prime_\infty$ denote the limit superconnection arising
as in the proof of Theorem \ref{th6}.
The spectral sequence for $\HH^*(A^\prime_\infty)$ gives
\begin{equation} \label{eq8.9}
\HH^1(A^\prime_\infty) \: = \: \HH^1(B; \R) \: \oplus \: \Ker \left( 
\HH^0(B; \HH^1(A^\prime_{\infty,[0]})) \rightarrow \HH^2(B; \R) \right).
\end{equation}
In particular,
\begin{equation}
\dim \left( \HH^1(A^\prime_\infty) \right) \: = \: \bb_1(B) \: + \: 
\dim \left( \Ker \left( 
\HH^0(B; \HH^1(A^\prime_{\infty,[0]})) \rightarrow \HH^2(B; \R) \right) 
\right).
\end{equation}
We wish to compare this with the corresponding spectral sequence for
$\HH^*(A^\prime_E)$, i.e. (\ref{eq8.11}).

Suppose that the differential $d_2 : \HH^0(B; \HH^1(Z; \R)) \rightarrow 
\HH^2(B; \R)$ vanishes. Then from (\ref{eq8.11}),
\begin{equation}
\bb_1(M) \: = \:
\bb_1(B) \: + \: \dim \left( \HH^0(B; \HH^1(Z; \R)) \right).
\end{equation}
By assumption, $\dim \left( \HH^1(A^\prime_\infty) \right) \: = \: j \: > \:
\bb_1(M)$. This implies that
\begin{equation} \label{ineq}
\dim \left( \HH^0(B; \HH^1(A^\prime_{\infty,[0]})) \right) \:
>  \: \dim \left( \HH^0(B; \HH^1(Z; \R)) \right).
\end{equation}

In terms of the original superconnection $A^\prime_E$,
we have that $\HH^1(Z; \R)$ is a flat subbundle of $\HH^1(A^\prime_{E,[0]})$. 
After taking limits, we obtain a flat subbundle 
$\HH^1(Z; \R)_\infty$ of $\HH^1(A^\prime_{\infty,[0]})$.
Here the fibers of $\HH^1(Z; \R)_\infty$ are again isomorphic to the
first real cohomology group of $Z$, but the flat structure could be different
than that of the bundle which we denoted by $\HH^1(Z; \R)$.
In particular,
\begin{equation}
\dim \left( \HH^0(B; \HH^1(Z; \R)_\infty) \right) \: \ge \: 
\dim \left( \HH^0(B; \HH^1(Z; \R)) \right).
\end{equation}
Clearly
\begin{equation}
\dim \left( \HH^0(B; \HH^1(A^\prime_{\infty,[0]})) \right) \:
\ge  \: \dim \left( \HH^0(B; \HH^1(Z; \R)_\infty) \right).
\end{equation}
Then from
(\ref{ineq}), we must have
\begin{equation} \label{ineq1}
\dim \left( \HH^0(B; \HH^1(Z; \R)_\infty) \right) \: > \: 
\dim \left( \HH^0(B; \HH^1(Z; \R)) \right).
\end{equation}
or 
\begin{equation} \label{ineq2}
\dim \left( \HH^0(B; \HH^1(A^\prime_{\infty,[0]})) \right) \:
>  \: \dim \left( \HH^0(B; \HH^1(Z; \R)_\infty) \right).
\end{equation}

If (\ref{ineq1}) holds then the
holonomy representation of the flat vector bundle $\HH^1(Z; \R)$ must have
a nontrivial unipotent subrepresentation (see 
\cite[Theorem 11.4 and
Proposition 11.14]{Richardson (1988)}). If (\ref{ineq2}) holds then
there is a nonzero covariantly-constant section of the vector bundle 
$\frac{H^1(A^\prime_{\infty,[0]})}{H^1(Z; \R)_\infty}$ on $B$, where the
flat connection on $\frac{H^1(A^\prime_{\infty,[0]})}{H^1(Z; \R)_\infty}$
is induced from the flat connection on 
$\HH^1(A^\prime_{\infty,[0]})$. This proves the corollary. $\square$ \\ \\
{\bf Proof of Corollary \ref{cor6} :}
Suppose that for $q \in \{p-1,p\}$, $b_q(Z) \: = \: \dim \left(
\Lambda^q({\frak n})^F \right)$.
From the Leray spectral
sequence, $\HH^p(M) \: \cong \: \Ker(\Phi^p - I) \oplus
\Coker( \Phi^{p-1} - I)$.
Let $\HH^*(Z; \R)_\infty$ denote the limiting
flat vector bundle on $S^1$, as in the proof of Corollary \ref{cor5},
with holonomy $\Phi^p_\infty \in
\Aut(\HH^p(Z; \R))$.
The spectral sequence for $\HH^*(A^\prime_\infty)$
gives $\HH^p(A^\prime_\infty) \: \cong \: \Ker(\Phi^p_\infty - I) \oplus
\Coker( \Phi^{p-1}_\infty - I)$. 
We have $\dim \left( \Ker(\Phi^p_\infty - I) \right) \ge 
\dim \left( \Ker(\Phi^p - I) \right)$ and
$\dim \left( \Coker(\Phi^p_\infty - I) \right) \ge 
\dim \left( \Coker(\Phi^p - I) \right)$.
By assumption, 
$j  =  \dim \left( \HH^p(A^\prime_\infty) \right) > \dim \left( 
\HH^p(M; \R) \right) = \bb_p(M)$. If
$\dim \left( \Ker(\Phi^p_\infty - I) \right) > 
\dim \left( \Ker(\Phi^p - I) \right)$ then $\Phi^p$ must have a
nontrivial unipotent subfactor.  Similarly, if 
$\dim \left( \Coker(\Phi^p_\infty - I) \right) > 
\dim \left( \Coker(\Phi^p - I) \right)$ then $\Phi^{p-1}$ must have a
nontrivial unipotent subfactor. $\square$ \\ \\
{\bf Example 9 : } Suppose that the affine fiber bundle
$M \rightarrow S^1$ has fiber $Z = T^2$. If $M$ has a $Sol$-geometry
or an $\R^3$-geometry
then Corollary \ref{cor6} implies that there are no small positive
eigenvalues in a
collapsing sequence associated to $M \rightarrow S^1$. On the other hand,
if $M$ has a $Nil$ geometry then Example 7 shows that there are small
positive
eigenvalues of the $1$-form Laplacian. See \cite{Jammes (2000)} for further
examples of homogeneous collapsings.\\ \\
{\bf Proof of Corollary \ref{cor7} :}
The $E_2$ term of the spectral sequence to compute $\HH^*(M; \R)$ consists
of $E_2^{p,0} = \HH^p(X; \R)$ and $E_2^{p,1} = \HH^p(X; {\cal O})$. The
differential is ${\cal M}_{\chi}$. The
corollary now follows from Theorem \ref{th6}. $\square$ \\ \\
{\bf Proof of Theorem \ref{th7} :}
As in \cite[\S 6]{Fukaya (1989)}, we can reduce the structure group of the 
fiber bundle $P \rightarrow B$ so that the local holonomy lies in a maximal
connected compact subgroup of $\Aff(Z)$, a torus group. Choose a horizontal
distribution $T^H M$ on $M$ whose local holonomy lies in this torus group. 
Add vertical Riemannian metrics $g^{TZ}$, parallel along the fibers, and a
Riemannian metric $g^{TB}$ on $B$ to give
$M \rightarrow B$ the structure of a Riemannian affine fiber bundle. We
will use Theorem \ref{th1} to make statements about the eigenvalues of the 
differential form Laplacian on $M$.

There is a vector space isomorphism 
${\frak n}^* \cong \bigoplus_{k=0}^S {\frak r}^*_{[k]}$. 
Define a number operator on ${\frak n}^*$ to be multiplication
by $3^k$ on ${\frak r}^*_{[k]}$. Extend this to a number operator on
$\Lambda^*({\frak n}^*)^F$ and to a number operator $N$ on the
vector bundle $E^*$ over $B$.

For $\epsilon > 0$, rescale $g^{TZ}$ to a new metric $g^{TZ}_\epsilon$
by multiplying it by 
$\epsilon^{3^k}$ on ${\frak r}_{[k]} \subset {\frak n}$. Let
$g^{TM}_\epsilon$ be the corresponding Riemannian metric on $M$.
The rescaling does not
affect $d^M$. The adjoint of $d^M$ with respect to the new metric
is $(d^M)^*_\epsilon \: = \: \epsilon^{N} \: (d^M)^* \: 
\epsilon^{-N}$.
Putting 
\begin{align} \label{eq8.16}
C^{\prime}_\epsilon & = \epsilon^{-N/2} \: d^M \: \epsilon^{N/2}, \\
C^{\prime \prime}_\epsilon 
& = \epsilon^{N/2} \: (d^M)^* \: \epsilon^{-N/2}, \notag
\end{align}
we have that $C^{\prime}_\epsilon$ is a flat degree-$1$ superconnection, 
with $C^{\prime \prime}_\epsilon$ being its adjoint with respect to 
$g^{TZ}$.
The Laplacian $\triangle^M$
coming from $g^{TM}_\epsilon$ is conjugate to
$C^{\prime}_\epsilon \: C^{\prime \prime}_\epsilon \: + \:
C^{\prime \prime}_\epsilon \: C^{\prime}_\epsilon$.

By \cite[\S 6]{Fukaya (1989)}, $\lim_{\epsilon \rightarrow 0}
(M, g^{TM}_\epsilon) = B$ with bounded sectional curvature in the limit. 
(The proof in \cite[\S 6]{Fukaya (1989)} uses a scaling
by $\epsilon^{2^k}$, but the proof goes through for a scaling by
$\epsilon^{3^k}$. The phrase ``The element $Y_i$ of ${\frak g}$, through
the right action of G, ...'' of
\cite[p. 349 b9]{Fukaya (1989)} should read ``... the left action of G, ...'') 
Let $A^\prime_\epsilon$ denote the superconnection on $E$
constructed by restricting
$C^\prime_\epsilon$ to the fiberwise-parallel forms.
We will show that $\lim_{\epsilon \rightarrow 0} A^\prime_\epsilon = 
\nabla^G$, the flat connection on $G$. The theorem will then follow from 
Theorem \ref{th1} and Lemma \ref{lemma?}.

Consider $A^\prime_{\epsilon,[0]}$.
It acts on a fiber of $E$ by $\epsilon^{-N/2} \: 
d^{\frak n} \: \epsilon^{N/2}$. Consider first its action on
a fiber $({\frak n}^*)^F \cong \bigoplus_{k=0}^S ({\frak r}^*_{[k]})^F$ of 
$E^1$. As
$d^{\frak n}$ acts on $\Lambda^1({\frak n}^*)$ by the dual of the Lie
bracket and $[{\frak r}_{[k]}, {\frak r}_{[l]}] \subset 
\bigoplus_{m > \max(k,l)}
{\frak r}_{[m]}$, we have $d^{\frak n} {\frak r}_{[m]}^* \subset 
\bigoplus_{k,l < m}
{\frak r}_{[k]}^* \wedge {\frak r}_{[l]}^*$. It follows that
\begin{equation} \label{eq8.17}
\epsilon^{-N/2} \: d^{\frak n}
 \: \epsilon^{N/2} \: : \: ({\frak r}_{[m]}^*)^F \rightarrow
({\frak r}_{[k]}^* \wedge {\frak r}_{[l]}^*)^F
\end{equation}
is $O \left( 
\epsilon^{(3^m - 3^k-3^l)/2} \right)$.
We obtain that the action of 
$A^\prime_{\epsilon,[0]}$ on $E^1$ is  $O \left( \epsilon^{1/2} \right)$
as $\epsilon \rightarrow 0$. A similar argument shows that the action of 
$A^\prime_{\epsilon,[0]}$ on $E^*$ is  $O \left( \epsilon^{1/2} \right)$.

Now consider $A^\prime_{\epsilon,[1]} \: = \: \epsilon^{-N/2} \: \nabla^E
 \: \epsilon^{N/2}$. Put $F_k^* = P \times_{Aff(Z)} ({\frak r}^*_{[k]})^F$, 
so that
$E^1 \cong \bigoplus_{k=0}^S F_k^*$. (Here the $*$ in $F_k^*$ denotes an
adjoint, not a $\Z$-grading.) Consider first the action of
$C^\prime_{\epsilon,[1]}$ on $C^\infty(B; E^1)$.
As the holonomy of $\nabla^E$ comes from an
$\Aff(Z)$ action, we have $\nabla^E  : C^\infty ( B; F_k^*) \rightarrow
\bigoplus_{l \le k}  C^\infty ( B; F_l^*)$ (see the proof of Proposition
\ref{prop3}). If $l < k$ then the component
\begin{equation} \label{eq8.18}
\epsilon^{-N/2} \: \nabla^E
 \: \epsilon^{N/2} \: : \: C^\infty ( B; F_k^*) \rightarrow
C^\infty ( B; F_l^*)
\end{equation}
of $\epsilon^{-N/2} \: \nabla^E \: \epsilon^{N/2}$
is $O \left( \epsilon^{1/2} \right)$.
On the other hand, the component 
$\nabla^E \: : \: C^\infty ( B; F_k^*) \rightarrow 
C^\infty ( B; F_k^*)$ is
the restriction of the flat connection $\nabla^G$ from $G^1$ to
$F_k^*$.  A similar argument applies to all of $E$ to show that
as $\epsilon \rightarrow 0$, $A^\prime_{\epsilon,[1]} = 
\nabla^G + O \left( \epsilon^{1/2} \right)$.

Finally, consider $A^\prime_{\epsilon,[2]} \: = \: \epsilon^{-N/2} \: i_T
 \: \epsilon^{N/2}$. The curvature $T$ of the fiber bundle $M \rightarrow B$
is independent of $\epsilon$. As $T$ acts by interior multiplication on
the fibers of $E$, the action of
$\epsilon^{-N/2} \: i_T \: \epsilon^{N/2}$ on 
$({\frak r}^*_{[k]})^F \subset E^1$ is $O \left( \epsilon^{3^k/2} \right)$. 
A similar argument applies to all of $E$ to show that
as $\epsilon \rightarrow 0$, $A^\prime_{\epsilon,[2]} = 
O \left( \epsilon^{1/2} \right)$.

The theorem follows. $\square$

\end{document}